\documentclass[preprint,12pt]{elsarticle}

\usepackage{graphicx}
\usepackage{amsmath,amssymb,amsthm}
\usepackage{subfig}
\usepackage{float}
\usepackage{subfloat}
\usepackage{mathptmx}
\usepackage{pgfplots}
\usepackage{multirow}
\usepackage{array}
\usepackage{caption}

\graphicspath{{./Figures/}}

\newtheorem{theorem}{Theorem}[section]
\newtheorem{definition}[theorem]{Definition}

\newtheorem{remark}[theorem]{Remark}

\newcommand{\set}[2]{\left\{{#1}\,:~{#2}\right\}}
\newcommand {\average}[1] {\mbox{$\left\{\!\!\left\{ #1 \right\}\!\!\right\}$}}
\newcommand {\jump}[1] {\mbox{$\left[\!\left[ #1 \right]\!\right]$}}

\newenvironment{defn*}{\begin{definition}}{\end{definition}}
\numberwithin{equation}{section}

\journal{}

\begin{document}

\begin{frontmatter}



\title{Optimal Control of Convective FitzHugh-Nagumo Equation}


\author[thk]{Murat Uzunca\corref{cor1}}
\ead{muzunca@thk.edu.tr}

\author[bal]{Tu\u{g}ba ~K\"u\c{c}\"ukseyhan}
\ead{guney.tugba@metu.edu.tr}

\author[iam]{Hamdullah  Y\"ucel}
\ead{yucelh@metu.edu.tr}

\author[iam,math]{B\"ulent Karas\"{o}zen}
\ead{bulent@metu.edu.tr}

\cortext[cor1]{Corresponding author}

\address[thk]{Department of Industrial Engineering, University of Turkish Aeronautical Association, Ankara, Turkey}
\address[iam]{Institute of Applied Mathematics, Middle East Technical University, 06800 Ankara, Turkey}
\address[bal]{Department of Mathematics,  Bal{\i}kesir  University,  10145 Bal{\i}kesir,  Turkey}

\begin{abstract}
We investigate smooth and sparse optimal control problems for convective FitzHugh-Nagumo equation with
travelling wave solutions in moving excitable media. The cost function includes distributed space-time and
terminal observations or targets. The state and adjoint equations are discretized in space by symmetric
interior point Galerkin (SIPG) method and by backward Euler method in time. Several numerical results are
presented for the control of the travelling waves. We also show numerically
the validity of the second order optimality conditions for the local solutions of the sparse optimal control
problem for vanishing Tikhonov regularization parameter.  Further, we  estimate the distance between the discrete control and associated local optima numerically by the help of the perturbation method and the smallest eigenvalue of the reduced Hessian.
\end{abstract}


\end{frontmatter}

\section{Introduction}
\label{introduction}

Spatially and temporally varying structures occur in form of Turing patterns,  traveling waves, fronts,  periodic pulses in many physical, chemical, and biological systems. They are described mathematically in form of coupled semi-linear partial differential equations (PDEs)  \cite{JDMurray_2002}. The  FitzHugh-Nagumo (FHN) equation is one of the most known generic  model in  physiology, describing complex wave phenomena
in excitable or oscillatory media. The most known type of
the FHN equation in the literature consists of a PDE with a non-monotone nonlinear term, known as activator equation, and  an ordinary differential equation (ODE), known as inhibitor equation. We call such kind of activator-inhibitor system as classical FHN equation. Another type of the activator-inhibitor system is diffusive FHN equation consisting of an activator PDE, and one or two inhibitor PDEs \cite{Karasozen15}. Recently, the convective FHN equation has been proposed as  a model for wave propagation in blood coagulation and bioreactor systems \cite{EAErmakova_MAPanteleev_EEShnol_2005,EAErmakova_EEShnol_MAPanteleev_AAButylin_VVolpert_FIAtaullakhanov_2009,AILonabov_TKStarozhilova_2005}. The presence of convective field leads to complex wave phenomena, like triggering and autonomous waves in a moving excitable media \cite{EAErmakova_EEShnol_MAPanteleev_AAButylin_VVolpert_FIAtaullakhanov_2009}.

%We here investigate a numerical treatment of optimal control problem governed by convective FitzHugh-Nagumo system. . In the objective functional  we add the $L^1$-norm of control as an extra term to exhibit the sparsity of control. As spatial dicretization, we use a discontinuous Galerkin discretization, while we use  standard backward Euler approach for temporal discretization. Numerical experiments are employed to demonstrate the limits of the model in terms of various perspectives.

%Spatially and/or temporally varying structures have been observed in variety of physical,   chemical  and biological systems operating far from equilibrium, for instance, pattern formation, blood coagulation \cite{FIAtaullakhanovVIZarnitsyna_2002}. Despite the great variability, such a kind of problems exhibit many common properties. Therefore, use of the simplest models still plays a vital role in the understanding of the mechanisms of propagation in excitable media. This class of models usually includes a nonlinear reaction between species or components and a diffusive transport of them due to density gradients. The simplest and the most widely used models of such reaction-diffusion systems is the model of FitzHugh-Nagumo \cite{JDMurray_2002}. However, some of these systems are described better by existence of flow and convection, for instance, blood coagulation or  . It is quite challengeable  to investigate the behaviour of such systems since

The classical FHN  equation and  its PDE part the Schl\"ogl or Nagumo equations were investigated theoretically and numerically  for the wave-type optimal control solutions \cite{RBuchholz_HEngel_EKammann_FTroltzsch_2013,ECasas_CRyll_FTroltzsch_2013,ECasas_CRyll_FTroltzsch_2015,CRyll_JLober_SMartens_HEngel_FTroltzsch_2016}.   Optimal control of semi-linear parabolic equation is an active research field with many applications in controlling pattern formation \cite{MStoll_JWPearson_PKMaini_2015} and feedback control of the monodomain equations in cardiac electrophysiology \cite{Breiten14} to give a few examples.

In this paper, we investigate the numerical treatment of optimal control problems governed by the convective FHN equation. The uncontrolled solutions of the convective FHN equation behave like travelling waves \cite{EAErmakova_EEShnol_MAPanteleev_AAButylin_VVolpert_FIAtaullakhanov_2009}. To control such a travelling waves, we use sparse control, which is a non-smooth $L^1$-control cost in addition to the $L^2$-control cost. When the control signals are localized in some regions of the space-time cylinder, sparse control
provides solutions without any a priori knowledge of these sub-areas. Sparse optimal control  was first investigated in \cite{GStadler_2009} for linear elliptic equations and later studied in \cite{ECasas_RHerzog_GWachsmuth_2012,ECasas_CRyll_FTroltzsch_2015,GWachsmuth_DWachsmuth_2010} for  semi-linear elliptic and parabolic equations.

%Controlling wave type solutions is numerically challenging in general. The results in \cite{RBuchholz_HEngel_EKammann_FTroltzsch_2013} show that some difficulties emerge in the presence of the unstable equilibrium point, even for 1D models.  When we add  an external influence such as convection,  it cause  qualitative changes in the behaviour of the system. Therefore, the choice of space discretization is quite important.

Here, we use symmetric interior penalty Galerkin (SIPG) method for space discretization. The discontinuous Galerkin (dG) methods are more stable for convection dominated problems than the continuous finite element methods and they do not require  the stabilization terms like the streamline upwind/Petrov-Galerkin method (SUPG). The dG methods have  several advantages compared to other numerical techniques such as finite volume and finite element methods; the trial and test spaces can be easily constructed, inhomogeneous boundary conditions and curved boundaries can be handled easily. They are also flexible in handling  non-matching grids and in designing hp-adaptive mesh refinement. The dG methods were successfully applied for the steady state \cite{DLeykekhman_MHeinkenschloss_2012a,HYucel_PBenner15,HYucel_BKarasozen_2014}, the time dependent linear convection-diffusion-reaction \cite{Akman14,TAkman_HYucel_BKarasozen_2014}, and  the semi-linear steady state  \cite{HYucel_PBenner_2015a,HYucel_MStoll_PBenner_2015b} optimal control problems. To solve the optimization problem, we here apply the \emph{optimize-discretize} approach. The first order optimality conditions are derived and then they are discretized by using the dG method. We employ the projected  nonlinear conjugate gradient (CG) method as an optimization algorithm \cite{WWHager_HZhang_2006}. We show the controllability of the traveling waves of the FHN equation with target functions in the whole space-time domain and at the final time with and without sparse controls. In addition, the performance of the method is  demonstrated  for convection dominated problems by increasing the wave velocity, for sparse controls with different sparse parameters.

Due to the semi-linearity of the FHN equation, the control problem is non-convex. Therefore, the fulfillment of the first necessary conditions does not imply the optimality. In order guarantee the optimality, the second-order sufficient optimality conditions (SSCs) have to be checked. In the recent years, the fulfillment of the SSCs for infinite dimensional and finite dimensional semi-linear PDE constrained optimal control problems has been  investigated extensively (see, e.g., \cite{ECasas_FTroltzsch_2015} for a recent overview). Except few examples with analytical solutions, it is not possible to prove the SSCs for the infinite dimensional problems since the unknown optimal solution is required. Therefore the finite dimensional approximations of the infinite dimensional problem are considered. Provided that the local minima satisfies the SSC, one can check the SSC numerically by finding a bound for the distance between the local minima and discrete solution \cite{ARoesch_DWachsmuth_2008a,ARoesch_DWachsmuth_2012}. For this purpose, the associated coercivity constant of the reduced Hessian operator is estimated numerically by computing its smallest eigenvalue. Similar techniques were applied to measure how far the control obtained by a reduced order optimization model is away from a local full control solution \cite{EKammann_FTroltzsch_SVolkwein,OLass_STrenz_SVolkwein_2015}. Moreover, the Tikhonov regularization parameter in the cost function  expresses the cost
of the control, and it increases the numerical stability of the optimal
solution. Recently the second optimality conditions have been investigated for semi-linear parabolic control problems with the objective function, not including the Tikhanov regularization term \cite{ECasas_FTroltzsch_2015}. We test the discrete optimization problem for vanishing Tikhonov regularization parameter as in \cite{ECasas_CRyll_FTroltzsch_2015}. The numerical results of the control of two dimensional waves confirm the convergence of the optimal solutions for vanishing Tikhonov regularization parameter as it was demonstrated for the one dimensional wave solutions of the classical FHN equation in \cite{ECasas_CRyll_FTroltzsch_2015,CRyll_JLober_SMartens_HEngel_FTroltzsch_2016}.

%We follow here the perturbation technique  for semi-linear elliptic control problems in  \cite{NArada_ECasas_FTroeltzsch_2002a}.

%The numerical implementation of such type problems in two-dimensional space is challengeable and requires a lot of %effort. To reducing running time without not losing so much in terms of accuracy, we apply a standard model %reduction  method, that is proper orthogonal decomposition (POD). Although we need a greater number of basis %functions comparing to the standard semi-linear parabolic control problems, our gains in terms of running time and %iteration number are promising.

The paper is organized as follows: In the next section the optimal control problem governed by the convective FHN equation is described as a model problem. We first prove the existence and uniqueness of a convective FHN equation, called as a state equation, by transforming into the one with monotone nonlinearity. Then, we introduce the first and second order optimality conditions. In Section~\ref{disc_state} we give a full discretization of the optimality system using the SIPG  method in space and the backward Euler discretization in time.
%Numerical treatment  of optimal control problem is given in section~\ref{disc_opg}.
In Section~\ref{numeric}, we discuss  some benchmark examples  with and without sparse controls. We investigate the effect of the Tikhonov parameter as it goes to zero. Further, with the help of the perturbation method and the smallest eigenvalue of the reduced Hessian, we find a bound for the distance between the local minima and discrete solution. The paper ends with some conclusions.

%%%%%%%%%%%%%%%%%%%%%%%%%%%%%%%%%%%%%%%%%%%%%%%%%%%%%%%%%%%%%%%%%%%%%%%%%%%%%%%%%%%%%%%%%%%%%%%%%%%%

\section{Optimal control of the convective FHN system}
\label{problem}

In this paper, we consider optimal control problems governed by the following convective FHN system:
\begin{subequations}\label{M1}
\begin{align}
 y_t(x,t) - d_y \Delta y(x,t) + \mathbf{V} \cdot \nabla y(x,t) + g(y(x,t)) + z(x,t)   &=  u(x,t)  &  \text{in   } Q_{T},  \\
z_t(x,t) - d_z \Delta z(x,t) + \mathbf{V} \cdot \nabla z(x,t) + \epsilon (z(x,t) - c_3 y(x,t))  &=  0  & \text{in   } Q_{T},
\end{align}
with homogeneous Neumann boundary conditions
\begin{equation}
\partial_{\emph{n}} y(x,t) = 0, \quad \partial_{\emph{n}} z(x,t)=  0   \qquad \qquad \quad   \text{on   } \Sigma_{T}^{N},
\end{equation}
Dirichlet boundary conditions
\begin{equation}
y(x,t) =  y_D(x,t),\quad  z(x,t)= z_D(x,t) \qquad  \text{on  } \Sigma_{T}^{D},
\end{equation}
and initial conditions
\begin{equation}
y(x,0) =  y_0(x), \quad \;  z(x,0) =z_0(x) \qquad \; \; \;  \text{in  } \Omega.
\end{equation}
\end{subequations}
In this setting, let $T>0$ be a fixed end time. We denote $Q_T$ the time space cylinder $Q_T:= \Omega \times (0,T)$, where $\Omega=(0,L) \times (0,H)$ is a bounded, Lipschitz domain. The lateral surface is denoted by
$\Sigma=\Gamma \times (0,T)$. We use the notation $\Sigma_T^D := \Gamma_D \times (0,T)$ and $\Sigma_T^N := \Gamma_N \times (0,T)$, where Dirichlet $\Gamma_D$ and Neumann $\Gamma_N$  boundaries, where the Dirichlet $y_D, z_D \in H^{3/2}(\Gamma_D)$ and the Neumann boundary conditions are prescribed. Moreover, the initial functions are given as $y_0, z_0 \in L^{\infty}(\Omega)$. We denote the outward unit normal vector and the associated outward normal derivative on $\partial \Omega$ by $\mathbf{n}$ and  $\partial_{\emph{n}}$, respectively. The diffusion coefficients are  denoted by $d_y$ and $d_z$. The parameters $c_3$ and $\epsilon$ are real constants. Further, the function $g(y)$ denotes the cubic polynomial nonlinearity
\begin{equation}\label{M2}
g(y)=c_{1} y(y-c_{2})(y-1)
\end{equation}
with the nonnegative real numbers $c_{i},  i=1,2$, which is monostable for $0 < c_1< 20$ and $c_2=0.02$  \cite{EAErmakova_EEShnol_MAPanteleev_AAButylin_VVolpert_FIAtaullakhanov_2009} in contrast to the bistable cubic nonlinearity for the Schl\"ogl equation \cite{RBuchholz_HEngel_EKammann_FTroltzsch_2013},  the classical FHN equation \cite{ECasas_CRyll_FTroltzsch_2013}, and the diffusive FHN equation \cite{Karasozen15}. The velocity field denoted by $\mathbf{V}= (V_{x_1}, V_{x_2})$ is given along the $x_{1}$-direction with a parabolic profile
\begin{equation}\label{M3}
V_{x_1}(x_2) = ax_{2}(H-x_{2}), \quad V_{max}=\frac{1}{4}aH^2, \quad a >0 , \quad V_{x_2} = 0,
\end{equation}
where $V_{\max}$ denotes the maximum wave speed of the velocity field. Moreover, the velocity field is the divergence free, i.e.,  $\mbox{div} \; {\mathbf V}=0$.

We also make the following assumption for the solutions $y$, $z$
\begin{align}\label{assnega}
0=:y_0 \leq y \leq y_1, \qquad 0=:z_0 \leq z \leq z_1  \quad \hbox{a.e.}\;\; Q_T,
\end{align}
which is admissible for the sake of physical realism. We note that the bounds are constant.

Here, we want to minimize an objective function of misfit type, i.e., the function is
designed to penalize deviations of the function values from the observed or measured data. We
formulate our minimization functional such as
\begin{align}\label{M4}
(P_{\mu}) \;\;  \left \{ \min \limits_{u \in \mathcal{U}_{ad}} J(u):= I(u) + \mu j(u), \right.
\end{align}
with
\begin{eqnarray*}
I(u) &=& \frac{\omega_Q^y}{2} \int \limits_{0}^{T} \int \limits_{\Omega} \big( y(x,t)-y_Q(x,t) \big)^2 \; dx\;dt
       + \frac{\omega_Q^z}{2} \int \limits_{0}^{T} \int \limits_{\Omega} \big( z(x,t)-z_Q(x,t) \big)^2 \; dx\;dt \\
     & &+ \frac{\omega_T^y}{2}  \int \limits_{\Omega} \big( y(x,T)-y_T(x,T) \big)^2 \; dx
       + \frac{\omega_T^z}{2}  \int \limits_{\Omega} \big( z(x,T)-z_T(x,T) \big)^2 \; dx  \\
     && + \frac{\omega_u}{2}   \int \limits_{0}^{T} \int \limits_{\Omega}  (u(x,t))^2 \; dx\;dt, \\
j(u) &=& \int \limits_{0}^{T} \int \limits_{\Omega} |u(x,t)| \; dx\;dt,
\end{eqnarray*}
where the pair $(y,z)$ denotes the solution of (\ref{M1}) associated to control $u$. In (\ref{M1}), the partial differential equation for $y$ is said  \emph{activator} equation, while the one for $z$  is called the \emph{inhibitor} equation. The functions $y_T, z_T \in L^{\infty}(\Omega)$ and $y_Q, z_Q \in L^{\infty}(Q_T)$ are the given targets or desired states. We have given constants $\omega_Q^y, \omega_Q^z, \omega_T^y, \omega_T^z$, Tikhonov regularization parameter $\omega_u \geq0$, and sparse parameter $\mu \geq 0$.\\

We consider the optimization problem (\ref{M4})  with pointwise box constraints
\begin{equation}\label{M5}
u \in \mathcal{U}_{ad} := \{ u \in L^{\infty}(Q_T): \; u_a \leq u(x,t) \leq u_b \;\; \hbox{for   a.e   } (x,t) \in Q_T \}
\end{equation}
with the real numbers $u_a \leq u_b$.

The  aim of the optimal control is to ensure that the state variables $y$ and $z$ are  as close as possible  to the desired or observed states by minimizing the
objective functional in the $L^2$ or $L^1$--norms. For well-defined optimal solutions, one has to show that there exists a unique solution $(y,z)$ of (\ref{M1}) for each $u \in \mathcal{U}_{ad}$. The FHN equation (\ref{M1}) belongs to the class of semi-linear parabolic equations with a non-monotone nonlinearity. The theory of the existence and uniqueness of  solutions $(y,z)$ of the state equation (\ref{M1}) is more delicate than the monotone nonlinearities \cite{RBuchholz_HEngel_EKammann_FTroltzsch_2013,ECasas_CRyll_FTroltzsch_2013}. Next section, we show the existence and uniqueness of the state equation (\ref{M1}) by transforming the state equation to one with the monotone nonlinearity for the convective FHN equation  \eqref{M1}.
\
%\begin{remark}
%To predict whether  given targets can be exactly reached or not is often difficult. When they are not reached,  the optimal control problem (\ref{M4}) will most likely be %unsolvable for $\omega_c=0$. Therefore, we need  the choice of $\omega_c >0$ with the box constraints. Then, optimal control problem is solvable, even if given targets can not %be reached exactly.
%\end{remark}

\subsection{Well-posedness of the state equation}
The existence and uniqueness of a weak solution for the Schl\"ogel and FHN equation was shown in \cite{RBuchholz_HEngel_EKammann_FTroltzsch_2013} by transforming   (\ref{M1}) into the one with monotone nonlinearity using the transformation $y=e^{\eta t} v$ with sufficiently large parameter $\eta$ \cite{RBuchholz_HEngel_EKammann_FTroltzsch_2013}. Here we apply the same technique for the FHN  equation with the convective term. Next, we construct the upper and lower solutions for the transformed equation, which yield pointwise bounds for the desired solution following \cite{CPao_1992}. These bounds are then used as an initial iterates to construct two monotonically convergent sequences. Finally, we show that their common limit is the unique solution of the transformed equation with the monotone nonlinearity.

Let us first perform the transformation of the state equation (\ref{M1}) by substituting $y=e^{\eta t} v$. Then, we obtain the following system:
\begin{subequations}\label{monotonestate}
\begin{eqnarray}
v_t - d_y \Delta v + \mathbf{V} \cdot \nabla v + e^{-\eta t} g(e^{\eta t} v) + \eta v  &=&  e^{-\eta t}(u-z)   \quad \text{in   } Q_{T},  \\
z_t - d_z \Delta z + \mathbf{V} \cdot \nabla z + \epsilon (z - c_3 e^{\eta t} v)  &=&  0  \qquad \qquad \quad \text{in   } Q_{T}
\end{eqnarray}
with the boundary conditions
\begin{equation}
\partial_{\emph{n}} v = 0, \quad \partial_{\emph{n}} z=  0 \quad  \text{on   } \Sigma_{T}^{N}, \quad \hbox{and} \quad
e^{\eta t} v =  y_D,\quad  z =  z_D  \quad  \text{on  } \Sigma_{T}^{D},\\
\end{equation}
and the initial conditions
\begin{equation}
v(x,0) =  y_0(x), \quad \;  z(x,0)= z_0(x) \qquad  \text{in  } \Omega.
\end{equation}
\end{subequations}

\indent Here, the nonlinear term
\[
\tilde{g} : v \longmapsto e^{-\eta t} g(e^{\eta t} v) + \eta v
\]
is a monotone non-decreasing function with respect to $v$ for all $(x,t)\in Q_T$. Moreover, it satisfies the following properties \cite[Sec.~4.3]{FTroeltzsch_2010a}:
\begin{description}
\item[(i)]  For every fixed $v \in \mathbb{R}$ is Lebesgue measurable in $Q_T$.
\item[(ii)] For almost all  $(x,t)\in Q_T$, $\tilde{g}$ is twice continuously differentiable with respect to $v$ and locally Lipschitz continuous of order 2 with respect to $v$, i.e., there exists $L(\rho)=6c_1 e^{2 \eta t} > 0 $ such that
\[
|\tilde{g}_{vv} (x,t,v_1) - \tilde{g}_{vv}(x,t,v_2)| \leq L(\rho) |v_1-v_2|
\]
holds with for all $v_1, v_2 \in \mathbb{R}$ with  $|v_i| \leq \rho$, $i=1,2$.
\end{description}
The nonlinearity is uniformly bounded and monotone increasing in the following sense:
\begin{description}
\item[(iii)] There exists a constant $C=c_1c_2+\eta+2 e^{\eta t}(c_1c_2+c_1)>0$ such that
\[
|\tilde{g}(x,t,0)| + |\tilde{g}_v(x,t,0)| + |\tilde{g}_{vv}(x,t,0)| \leq C.
\]
\item[(iv)] It holds $0 \leq \tilde{g}_v (x,t,v)$  for almost all $(x,t)\in Q_T$, all $ v \in \mathbb{R}$.
\end{description}

Before defining a weak formulation of the system (\ref{monotonestate}), we need to define the following Hilbert space
\[
W(0,T):=\{w \in L^2(0,T; V);\; w' \in L^2(0,T; V^{*})\}
\]
equipped with the norm
\[
\|w\|_{W(0,T)}= \left( \int_0^T (\|w(t)\|^2_V + \|v_t(t)\|^2_{V^*})\mathrm{d}t \right)^{\frac{1}{2}},
\]
where $V=H^1(\Omega)$ and $V^*$ is the dual space of $V$. Now, we can define a weak solution of the system \eqref{monotonestate}.

\begin{definition}
A pair of functions $(v,z)\in (W(0,T) \cap L^{\infty}(Q_T))^2$ is called weak solution of the system \eqref{monotonestate}, if the equations
\begin{subequations}\label{weaksoln}
\begin{align}
\hspace{-4.5mm}\int_0^T (v_t, \varphi)_{V^*,V} \; \mathrm{d}t + \iint_{Q_T} \big( d_y \nabla v \cdot \nabla \varphi + \mathbf{V} \cdot \nabla v \varphi + \tilde{g} \varphi - e^{\eta t}(u-z) \varphi \big) \; \mathrm{d}x \;\mathrm{d}t =& 0, \\
\hspace{-4.5mm}\int_0^T (z_t, \varphi)_{V^*,V} \; \mathrm{d}t + \iint_{Q_T} \big( d_z \nabla z \cdot \nabla \varphi + \mathbf{V} \cdot \nabla z  \varphi + \epsilon (z-c_3 e^{\eta t}v) \varphi \big) \; \mathrm{d}x \; \mathrm{d}t =& 0,
\end{align}
and
\begin{equation}
v(\cdot, 0)=y_0, \qquad  z(\cdot, 0)=z_0
\end{equation}
\end{subequations}
are satisfied for all $\varphi \in L^2(0,T; H^1(\Omega))$. It is noted that $\nabla$ denotes the gradient with respect to $x$.
\end{definition}

Next, we give the definition of the ordered upper and lower solutions as done in \cite{WBarthel_CJohn_FTroltzsch_2010,CPao_1992}. The pair functions  $(\widetilde{v}, \widetilde{z})$ and  $(\widehat{v}, \widehat{z})$ in $(W(0,T) \cap L^{\infty}(Q_T))^2$ are said to be  ordered upper and lower solutions of \eqref{monotonestate}, respectively, if they satisfy
\[
(\widehat{v}, \widehat{z})  \leq (\widetilde{v}, \widetilde{z})
\] and
\begin{eqnarray*}
\widehat{v}_t - d_y \Delta \widehat{v} + \mathbf{V} \cdot \nabla \widehat{v} + \widetilde{g}(x,t,\widehat{v})- e^{-\eta t}(u-\widetilde{z}) \leq \; & 0 & \; \leq \widetilde{v}_t - d_y \Delta \widetilde{v} + \mathbf{V} \cdot \nabla \widetilde{v} \\
&& \quad + \widetilde{g}(x,t,\widetilde{v}) - e^{-\eta t}(u-\widehat{z}),\\
\widehat{z}_t - d_z \Delta \widehat{z} + \mathbf{V} \cdot \nabla \widehat{z} + \epsilon (z-c_3 e^{\eta t} \widetilde{v}) \leq \;& 0 &\; \leq \widetilde{z}_t - d_z \Delta \widetilde{z} + \mathbf{V} \cdot \nabla \widetilde{z} \\
&& \quad + \epsilon (z-c_3 e^{\eta t} \widehat{v}), \\
\partial_{\emph{n}} \widehat{v} \leq \;& 0 &\; \leq \partial_{\emph{n}} \widetilde{v}, \\
\partial_{\emph{n}} \widehat{z} \leq \;& 0 &\; \leq \partial_{\emph{n}} \widetilde{z}, \\
\widehat{v} \leq \;& e^{-\eta t} y_{D} &\; \leq \widetilde{v}, \\
\widehat{z} \leq \;& z_{D} &\; \leq \widetilde{z},\\
\widehat{v}(x,0) \leq \;& y_0(x) &\; \leq \widetilde{v}(x,0), \\
\widehat{z}(x,0) \leq \;& z_0(x) &\; \leq \widetilde{z}(x,0).
\end{eqnarray*}
By taking
\[
\widetilde{v}(x,t)=\widetilde{z}(x,t)=M, \qquad \widehat{v}(x,t)=\widehat{z}(x,t)=0
\]
for some $M>0$, we can rewrite our system \eqref{monotonestate} as
\begin{subequations}\label{remonotonestate}
\begin{eqnarray}
v_t - d_y \Delta v + \mathbf{V} \cdot \nabla v + \tilde{g}(e^{\eta t v}) - e^{-\eta t}(u-M) &=& e^{-\eta t}(M - z), \\
z_t - d_z \Delta z + \mathbf{V} \cdot \nabla z + \epsilon (z - c_3 e^{\eta t}M) &=& \epsilon c_3 e^{\eta t} (M-v), \\
\partial_{\emph{n}} v = 0, \quad \partial_{\emph{n}} z &=&  0, \\
v= e^{-\eta t} y_D, \quad z &=&z_D, \\
v(x,0) =  y_0(x), \quad \;  z(x,0) &=& z_0(x).
\end{eqnarray}
\end{subequations}
Here, we have
\begin{eqnarray*}
\frac{\partial (e^{-\eta t}(M - z)) }{\partial v}, \quad \frac{\partial (\epsilon c_3 e^{\eta t} (M-v)) }{\partial z} &\geq & 0, \\
\frac{\partial (e^{-\eta t}(M - z)) }{\partial z}, \quad \frac{\partial (\epsilon c_3 e^{\eta t} (M-v)) }{\partial v} &\leq & 0
\end{eqnarray*}
for all $v,z \in [0,M]$.

Now, we can state the existence and uniqueness of the system \eqref{monotonestate} for each control variable $u$.

\begin{theorem}
Assume that the initial conditions $y_0$ and $z_0$ are nonnegative functions, and (\ref{assnega}) holds. Then, the system \eqref{monotonestate} admits a unique solution $(v,z) \in (W(0,T)\cap C(Q_T))^2$ for each control $u \in \mathcal{U}_{ad}$.
\end{theorem}
\begin{proof}
We adopt the iteration technique introduced in \cite{CPao_1992} and construct sequences $\{ (\bar{v}\;^k, \bar{z}\;^k)\}_{k=0}^{\infty}$, $\{ (\underbar{v}\;^k, \underbar{z}\;^k) \}_{k=0}^{\infty}$ with initial elements
\begin{align*}
& \bar{v}\;^0 = \widetilde{v} = M, \qquad  \bar{z}\;^0 = \widetilde{z} = M, \\
& \underbar{v\;}^0=\widehat{v}=0,  \qquad    \underbar{z}\;^0=\widehat{z}=0.
\end{align*}
Initiating from $(\bar{v}\;^{k}, \bar{z}\;^{k})$ and $(\underbar{v}\;^{k}, \underbar{z}\;^{k})$, $(\bar{v}\;^{k+1}, \bar{z}\;^{k+1})$ and $(\underbar{v}\;^{k+1}, \underbar{z}\;^{k+1})$ are found by solving
\begin{eqnarray*}
\bar{v}\;_t^{k+1} - d_y \Delta \bar{v}\;^{k+1} + \mathbf{V} \cdot \nabla \bar{v}\;^{k+1} + \tilde{g}(e^{\eta t} \bar{v}\;^{k+1}) - e^{-\eta t}(u-M) &=& e^{-\eta t}(M - \underbar{z}\;^k), \\
\bar{z}\;_t^{k+1} - d_z \Delta \bar{z}\;^{k+1} + \mathbf{V} \cdot \nabla \bar{z}\;^{k+1} + \epsilon (\bar{z}\;^{k+1} - c_3 e^{\eta t}M) &=& \epsilon c_3 e^{\eta t} (M-\underbar{v}\;^k), \\
\partial_{\emph{n}} \bar{v}\;^{k+1} = 0, \quad \partial_{\emph{n}} \bar{z}\;^{k+1} &=&  0, \\
\bar{v}\;^{k+1} = e^{-\eta t} y_D, \quad \bar{z}\;^{k+1} &=& z_D, \\
\bar{v}\;^{k+1}(x,0) =  y_0(x), \quad \;  \bar{z}\;^{k+1}(x,0) &=& z_0(x)
\end{eqnarray*}
and
\begin{eqnarray*}
\underbar{v}\;_t^{k+1} - d_y \Delta \underbar{v}\;^{k+1} + \mathbf{V} \cdot \nabla \underbar{v}^{k+1} + \tilde{g}(e^{\eta t} \underbar{v}\;^{k+1}) - e^{-\eta t}(u-M) &=& e^{-\eta t}(M - \bar{z}\;^k),  \\
\underbar{z}\;_t^{k+1} - d_z \Delta \underbar{z}\;^{k+1} + \mathbf{V} \cdot \nabla \underbar{z}\;^{k+1} + \epsilon (\underbar{z}\;^{k+1} - c_3 e^{\eta t}M) &=& \epsilon c_3 e^{\eta t} (M-\bar{v}\;^k), \\
\partial_{\emph{n}} \underbar{v}\;^{k+1} = 0, \quad \partial_{\emph{n}} \underbar{z}\;^{k+1} &=&  0, \\
\underbar{v}\;^{k+1} = e^{-\eta t} y_D, \quad \underbar{z}\;^{k+1} &=& z_D, \\
\underbar{v}\;^{k+1}(x,0) =  y_0(x), \quad \;  \underbar{z}\;^{k+1}(x,0) &=& z_0(x),
\end{eqnarray*}
respectively. The constructed sequence  $\{ (\bar{v}\;^k, \bar{z}\;^k)\}_{k=0}^{\infty}$ is monotone non-increasing and upper solution for all $k$. Conversely $\{ (\underbar{v}\;^k, \underbar{z}\;^k) \}_{k=0}^{\infty}$ is monotone non-decreasing and lower solution for all $k$. Further, we have
\[
\underbar{u}\;^k (x,t) \leq \bar{u}\;^k (x,t) \quad \text{and} \quad \underbar{v}\;^k (x,t) \leq \bar{v}\;^k (x,t)
\]
for all  $k \in \mathbb{N}$ and $(x,t)\in Q_T$.

By induction, we can verify the monotonicity of the sequence $\{\bar{v}\;^k\}_{k=0}^{\infty}$: For $k=0$,
\begin{eqnarray*}
\bar{v}\;_t^1 - d_y \Delta \bar{v}\;^1 + \mathbf{V} \cdot \nabla \bar{v}\;^1 + \tilde{g}(e^{\eta t} \bar{v}\;^{1}) - e^{-\eta t}(u-M) &=& e^{-\eta t}(M - \underbar{z}\;^0),  \\
\bar{z}\;_t^1 - d_z \Delta \bar{z}\;^1 + \mathbf{V} \cdot \nabla \bar{z}\;^1 + \epsilon (\bar{z}\;^1 - c_3 e^{\eta t}M) &=& \epsilon c_3 e^{\eta t} (M-\underbar{v}\;^0), \\
\partial_{\emph{n}} \bar{v}\;^1 = 0, \quad \partial_{\emph{n}} \bar{z}\;^1 &=& 0, \\
\bar{v}\;^1 = e^{-\eta t} y_D, \quad \bar{z}\;^1 &=& z_D, \\
\bar{v}\;^1(x,0) =  y_0(x), \quad \;  \bar{z}\;^1(x,0) &=& z_0(x).
\end{eqnarray*}
The property that $\bar{v}\;^0$ is an upper solution  gives us
\begin{eqnarray*}
\bar{v}\;_t^0 - d_y \Delta \bar{v}\;^0 + \mathbf{V} \cdot \nabla \bar{v}\;^0 + \tilde{g}(e^{\eta t} \bar{v}\;^{0}) - e^{-\eta t}(u-M) &\geq& e^{-\eta t}(M - \underbar{z}\;^0),  \\
\bar{z}\;_t^0 - d_z \Delta \bar{z}\;^0 + \mathbf{V} \cdot \nabla \bar{z}\;^0 + \epsilon (\bar{z}\;^0 - c_3 e^{\eta t}M) &\geq& \epsilon c_3 e^{\eta t} (M-\underbar{v}\;^0), \\
\partial_{\emph{n}} \bar{v}\;^0 \geq 0, \quad \partial_{\emph{n}} \bar{z}\;^0 &\geq&  0, \\
\bar{v}\;^0 \geq  e^{-\eta t} y_D, \quad \bar{z}\;^0 & \geq & z_D, \\
\bar{v}\;^0(x,0) \geq  y_0(x), \quad \;  \bar{z}\;^0(x,0) &\geq& z_0(x).
\end{eqnarray*}
Hence, we obtain
\begin{eqnarray*}
\bar{v}\;_t^0 -\bar{v}\;_t^1  - d_y \Delta (\bar{v}\;^0-\bar{v}\;^1) + \mathbf{V} \cdot \nabla (\bar{v}\;^0-\bar{v}\;^1) + \tilde{g}(e^{\eta t} \bar{v}\;^{0}) - \tilde{g}(e^{\eta t} \bar{v}\;^{1}) &\geq& 0,  \\
\bar{z}\;_t^0-\bar{z}\;_t^1 - d_z \Delta (\bar{z}\;^1-\bar{z}\;^0) + \mathbf{V} \cdot \nabla (\bar{z}\;^0-\bar{z}\;^1) + \epsilon (\bar{z}\;^0 - \bar{z}\;^1) &\geq& 0, \\
\partial_{\emph{n}} (\bar{v}^0-\bar{v}^1) \geq 0, \quad \partial_{\emph{n}} (\bar{z}^0-\bar{z}^1) &\geq&  0, \\
(\bar{v}\;^0-\bar{v}\;^1) \geq 0, \quad (\bar{z}\;^0-\bar{z}\;^1) &\geq& 0, \\
(\bar{v}\;^0-\bar{v}\;^1)(x,0) \geq  0, \quad \;  (\bar{z}\;^0-\bar{z}\;^1)(x,0) &\geq& 0.
\end{eqnarray*}
So it follows from the comparison principle for nonlinear parabolic equations $\bar{v}\;^0 - \bar{v}\;^1,\; \bar{z}\;^0 - \bar{z}\;^1 \geq 0$. Now if we assume that $\bar{v}\;^{k-1} - \bar{v}\;^k, \; \bar{z}\;^{k-1} - \bar{z}\;^k \geq 0$, one can easily show that $\bar{v}\;^k - \bar{v}\;^{k+1}, \; \bar{z}\;^k - \bar{z}\;^{k+1} \geq 0$. Analogously, the monocity of $(\underbar{v}\;^k,\underbar{z}\;^k )$ can be proved.

Now, we show the convergence of  the sequence $\{\bar{v}\;^k, \bar{z}\;^k \}$ to a solution of (\ref{monotonestate}). The sequence $\{\bar{v}\;^k\}$ is monotone non-increasing and bounded from below by $\hat{u}=0$. Hence by Lebesgue dominated convergence theorem \cite{WBarthel_CJohn_FTroltzsch_2010,CPao_1992}, it converges to $v$ in space $L^p(Q)$, $p < \infty$. On the other hand, the sequence $\{\bar{z}\;^k\}$ is monotone non-decreasing and bounded from above by $\tilde{z}=M$. It converges to $z$ in a similar way.

Finally, we prove the uniqueness of the solution of (\ref{monotonestate}). Suppose that $(v_1, z_1)$, $(v_2, z_2)$ are solutions of \eqref{weaksoln} and set $v:=v_1-v_2,\; z:=z_1-z_2$.  Then $v$ and $z$ satisfy the initial conditions obviously. Moreover, the following equations
\begin{subequations}\label{weakfor uniq}
\begin{align}
\int_0^T (v_t, \varphi)_{V^*,V} \;\mathrm{d}t +& \iint_{Q_T} (d_y \nabla v \cdot \nabla \varphi + \mathbf{V} \cdot \nabla v  \varphi + \tilde{g}(x,t,v_1) \varphi \nonumber \\
  &- \tilde{g}(x,t,v_2) \varphi - e^{\eta t}(u-z) \varphi) \; \mathrm{d}x \; \mathrm{d}t = 0, \label{weakfor uniqa} \\
\int_0^T (z_t, \varphi)_{V^*,V} \; \mathrm{d}t +& \iint_{Q_T} (d_z \nabla z \cdot \nabla \varphi + \mathbf{V} \cdot \nabla z  \varphi \nonumber \\
&+ \epsilon (z-c_3 e^{\eta t}v) \varphi) \; \mathrm{d}x \; \mathrm{d}t = 0 \label{weakfor uniqb}
\end{align}
\end{subequations}
hold for all $\varphi\in W(0,T))$. Then, following the steps in \cite{RGriesse_2003} by taking
$\varphi=v$ in (\ref{weakfor uniqa}) and  $\varphi=z$ in  (\ref{weakfor uniqb}) we obtain the desired result $v=0$ and $z=0$.
\end{proof}

Hence, we can give the existence of an optimal control $u$ for the optimal control problem (\ref{M4}).

\begin{theorem}
The optimal control problem (\ref{M4}) has at least one optimal solution $u$  with associated optimal state $y$.
\end{theorem}
\begin{proof}
Here we just sketch the key ideas of the proofs in \cite[Sec.~5.3,Thm.~7.4]{FTroeltzsch_2010a,RGriesse_2003}. Since $\mathcal{U}_{ad}$ is non-empty and bounded in $L^{\infty}(Q_T)$, it is bounded in any space $L^{p}(Q_T)$ and it follows from the existence and uniqueness of the state variables that they are also  bounded. Hence, the cost functional is bounded below, which allows the existence of an infimum. Therefore one can find a weakly convergent minimizing sequence due to the boundedness of this sequence. Then, the compact embedding results give us the strong convergence of the state in a weaker norm. Hence, there exists a feasible limit point, and convergence of the objective function can be shown using the continuity argument.
\end{proof}

\begin{remark}
In this paper, we show the well-posedness of the state equation  by introducing upper and lower solutions for the state with monotone nonlinearity, obtained after a suitable transformation. However, there exists other possible reformulations in the literature, for instance, Schauder fixed point theorem applied in \cite{DEJackson_1990} for FitzHugh-Nagumo equation, Brouwer fixed point theorem applied in \cite{TGudi_AKPani_2007} for a class of quasi--linear elliptic problems which are of nonmonotone type, Leray-Schauder fixed-point theorem applied in \cite{RGriesse_SVolkwein_2005a} for a coupled system of semi--linear parabolic reaction--diffusion equations, Faedo-Galerkin method applied in \cite{JYPArk_SHPark_2007} for  a hyperbolic quasi-linear hemivariational inequalities.
\end{remark}

We continue this section by introducing the necessary and sufficient optimality conditions of the optimal control problem (\ref{M4}).

\subsection{First order optimality conditions}
The OCP
\begin{align}\label{M6}
 \min \limits_{u \in \mathcal{U}_{ad}} J(u):=f(y_u,z_u,u)= I(u) + \mu j(u)
\end{align}
subject to the convective FHN equation (\ref{M1}) is a non-convex programming problem so that  different local minima might occur. Since any global solution is one of these local solutions, we set up the first-order optimality conditions satisfied by the local minima.

The cost functional $J(u)$ in (\ref{M6}) consists of two terms with different smoothness. While the first part $I(u)$ is smooth, the second part $j(u): L^1(Q) \rightarrow \mathbb{R}$ is not  differentiable, but it is subdifferentiable and  the directional derivative is given by
\begin{equation}\label{M7}
j^{\;'}(u, v- u) = \max \limits_{\lambda \in \partial j(u)} < \lambda, v-u>
\end{equation}
with
\[
\partial j(u) = \left \{ \lambda \in L^{\infty}(Q_T)\; : j(v) \geq j(u) + \int \limits_{0}^{T} \int \limits_{\Omega} \lambda (v- u) \; dx\; dt  \quad \forall v \in L^{\infty}(Q_T) \right \},
\]
where
\[
\lambda(x,t) \in \left\{
  \begin{array}{ll}
    \{1\},  & \hbox{if} \;\; v(x,t) >0, \\
    \left [-1. 1\right ], & \hbox{if} \;\; v(x,t) =0, \\
    \{-1\}, & \hbox{if} \;\; v(x,t) <0.
  \end{array}
\right.
\]
We note that the relations above are required only almost everywhere.  With the help of  the Lagrangian functional, we obtain the following variational inequality:
\begin{align}\label{M8}
I^{\;'}(u)(v-u) + \mu j^{\;'}(u,v-u) \geq 0  \qquad \forall v \in \mathcal{U}_{ad},
\end{align}
that is,
\begin{align*}
\int \limits_{0}^{T} \int \limits_{\Omega} \big( p(x,t) + \omega_u u(x,t) + \mu \lambda(x,t) \big) \big(v(x,t) -u(x,t)\big) \;dx \;dt \geq 0 \qquad \forall v \in \mathcal{U}_{ad},
\end{align*}
where  $p(x,t)$ with $q(x,t)$ are called the adjoint variables as the solution of the following system
\begin{subequations}\label{M9}
\begin{eqnarray}
-p_t - d_y\Delta p - \mathbf{V} \cdot \nabla p + g_y(y) p - \epsilon c_3 q  &=& \omega_Q^y \big(y - y_Q \big) \quad  \text{in  } Q_T,   \\
-q_t - d_z\Delta q - \mathbf{V} \cdot \nabla q + \epsilon q + p  &=&  \omega_Q^z \big(z-z_Q \big) \quad   \text{in  } Q_T,
\end{eqnarray}
with the mixed boundary conditions
\begin{equation}
d_y \partial_{\emph{n}} p(x,t) + (\mathbf{V}\cdot n) p(x,t) = 0, \quad  d_z \partial_{\emph{n}} q(x,t) + (\mathbf{V} \cdot n)q(x,t)  =  0 \quad \text{on  } \Sigma_T^N,
\end{equation}
the Dirichlet boundary conditions
\begin{equation}
p(x,t) =  0,\qquad  q(x,t)  =  0 \qquad  \text{on  } \Sigma_T^D,  \\
\end{equation}
and final time conditions
\begin{equation}
p(x,T) =  \omega_T^y \big( y(x, T)-y_T(x) \big), \quad   q(x,T)  = \omega_T^z \big( z(x, T)-z_T(x) \big) \quad  \text{in  } \Omega. \end{equation}
\end{subequations}
The convection term in the adjoint system (\ref{M9}) is the negative of the one in the FHN equation (\ref{M1}). As a consequence, errors in the solution can potentially propagate in both directions. Therefore, the numerical treatment of the state and adjoint systems together is more delicate.

For $\omega_u >0$ and  $\mu >0$, from the variational inequality (\ref{M8}) the following  projection formulas  are obtained \cite{ECasas_CRyll_FTroltzsch_2015,CRyll_JLober_SMartens_HEngel_FTroltzsch_2016}
\begin{eqnarray}
u(x,t) &=& \mathbb{P}_{[u_a,u_b]} \left\{ -\frac{1}{\omega_u} \big(p(x,t) + \mu \lambda(x,t) \big) \right\} \quad \hbox{for  a.a  } (x,t) \in Q_T , \label{M10} \\
\lambda(x,t) &=& \mathbb{P}_{[-1,1]} \left\{ -\frac{1}{\mu} p(x,t) \right\}  \quad \hbox{for  a.a  }   (x,t) \in Q_T, \label{M11}
\end{eqnarray}
where the projection operator $\mathbb{P}_{[a,b]}: \mathbb{R} \rightarrow [a,b]$ is defined by
\[
\mathbb{P}_{[a,b]}(f)=\max \{a, \min\{f,b\} \}.
\]
Further, the following relation holds for  almost all $(x,t) \in Q_T$
\begin{equation}\label{M12}
u(x,t) = 0 \;\; \Leftrightarrow \;\; \left\{
                                       \begin{array}{ll}
                                         |p(x,t)| \leq \mu, & \hbox{if} \;\; u_a <0, \\
                                         p(x,t) \geq - \mu, & \hbox{if} \;\; u_a =0.
                                       \end{array}
                                     \right.
\end{equation}
We refer to \cite{ECasas_RHerzog_GWachsmuth_2012,ECasas_CRyll_FTroltzsch_2013,FTroeltzsch_2010a} for a further discussion on the projection operator and a derivation of (\ref{M11}) and (\ref{M12}).

\subsection{Second order optimality conditions}

%We are interested  in the numerical solution of the optimal control problems governed by the convective FHN equation (\ref{M4}).
 Due to the nonlinearity of the state equation (\ref{M1}), the optimization problem is non-convex. Therefore, the fulfillment of the first necessary conditions does not imply optimality. In order guarantee the optimality, the second-order sufficient optimality conditions (SSCs) have to be satisfied. The SSCs are related to certain critical cones that must be chosen as small as possible.

%Here we investigate SSCs by measuring the distance of discrete solution $u_h$ of (\ref{M4}) near a local solution %$u$ of (\ref{M4}).

Now, we give the critical cone related to our optimization problem  for $\omega_u > 0 $ derived in \cite{ECasas_CRyll_FTroltzsch_2015}:
\[
C_u = \left \{ v \in L^2(Q_T): \; v \hbox{  satisfies the sign condition   and   } I^{'}(u)v + \mu j^{\,'}(u)=0 \right \}
\]
with  the sign condition:
\[
\left\{
  \begin{array}{ll}
    v(x,t) \geq 0, & \hbox{if} \;\; u(x,t)=u_a, \\
    v(x,t) \leq 0, & \hbox{if} \;\; u(x,t)=u_b.
  \end{array}
\right.
\]
The set $C_u$ is a convex and closed cone in $L^2(Q_T)$. Moreover, if $u$ is a local minima for $(P_{\mu})$, then the following inequalities hold \cite[Theorem 3.3]{ECasas_CRyll_FTroltzsch_2015},
\begin{equation}\label{ssc}
I^{\,''}(u) v^2 \geq 0 \quad \forall v \in C_u \backslash \{0\},  \;\; \hbox{equivalently} \;\;
I^{\,''}(u) v^2 \geq \delta \|v\|^2_{L^2(Q_T)} \;\;  \forall v \in C_u, \;\; \delta>0.
\end{equation}
Then, under the assumption $I^{\,''}(u) v^2 \geq 0 \quad \forall v \in C_u \backslash \{0\}$, there exits $\delta>0$ and $r_0>0$ such that
\begin{equation}\label{ssc1}
J(u) + \frac{\delta}{2} \|v-u\|^2_{L^2(Q_T)} \leq J(v) \quad \forall u \in U_{ad} \cap B_{r_0}(u),
\end{equation}
where $B_{r_0}(u)$ is the $L^2(Q_T)$ ball centered at $u$ with radius $r_0$. This shows the existence of a local minima, see \cite[Theorem 3.4]{ECasas_CRyll_FTroltzsch_2015} for details.

%We need to compute the second-order condition (\ref{ssc}) to measure the distance a numerical solution $u_h$ of a dicretization (\ref{M4}) and a local minima $u$.
The verification of the SSCs is difficult because the solution of the infinite dimensional problem is required. Even if it is known, it would still be tedious to check that the SSC holds, since it requires the exact solution of linearized PDEs. However, there exist some numerical studies on the SSCs, see e.g., \cite{EKammann_FTroltzsch_SVolkwein,HDMittelmann,ARoesch_DWachsmuth_2008a}. Here we determine the constant $\delta$ by computing the smallest eigenvalue of the reduced Hessian as introduced in \cite{HDMittelmann} of the finite dimensional dG discretized OCP.

%for some OCPs with given analytical solutions it can be verified

%In general, we are aware of that this way is not reliable in estimating the constant $\delta$ for the infinite dimensional optimal control problem. However, there is no reliable and %practical method to verify the condition (\ref{ssc}).

Now, we can state the following theorem to measure the distance between the local minima $u$ and the local discrete solution $u_h$ obtained by applying the discontinuous Galerkin discretization for spatial discretization  and the  backward Euler for temporal discretization in Section~\ref{disc_state}.

\begin{theorem}\label{Thm:est}
Let $u$ be a local minima of (\ref{M4}). Assume that $u$ satisfies the second-order sufficient condition (\ref{ssc}) and (\ref{ssc1}). If $u_h$ is the discrete solution such that $u_h \in B_{r_0}(u)$, then it holds
\begin{equation}
\|u - u_h\|_{L^2(Q_T)} \leq \frac{1}{\delta^{\;'}} \|\zeta\|_{L^2(Q_T)},
\end{equation}
where a perturbation function $\zeta$ is defined as the following
\[
\zeta(x):=\left\{
  \begin{array}{ll}
    -\min\{0,\omega_u u_h + p_h + \mu \lambda_h\}, & \hbox{if}  \;\;  u_{h} = u_a, \\
    - \big( \omega_u u_h + p_h + \mu \lambda_h \big), & \hbox{if}  \;\; u_a < u_h < u_b, \\
    -\max\{0,\omega_u u_h + p_h + \mu \lambda_h\}, & \hbox{if}  \;\; u_h = u_b
  \end{array}
\right.
\]
and  $0 < \delta^{\;'} <\delta$.
\end{theorem}
\begin{proof}
Let $u_h$ be a discrete solution that need not be optimal for the continuous problem $(\ref{M4})$, and let $p_h$ and $\lambda_h$ be the associated adjoint and sparse. If $u_h$ were optimal, then $p_h + \omega_u u_h + \mu \lambda_h =0$ should be satisfied in  almost all points  $ x \in \Omega$, where $u_a \leq u_h \leq u_b$ holds. If not, then $p_h + \omega_u u_h + \mu \lambda_h + \zeta =0$, where $\zeta$ is a perturbation function, adopted from \cite{NArada_ECasas_FTroeltzsch_2002a}. Although $u_h$ is possibly not be optimal for (\ref{M4}), it is optimal for the perturbed optimization problem
\[
\min \limits_{u \in \mathcal{U}_{ad}} \; J(u) + (\zeta, u)_{L^2(Q_T)}.
\]
Inserting $u$ in the discrete variational form and $u_h$ in the continuous discrete variational form,  we obtain
\begin{subequations}\label{sm}
\begin{eqnarray}
\left( J^{\;'}(u_h) + \zeta, u-u_h \right) &\geq& 0, \\
\left( J^{\;'}(u) + \zeta, u_h-u \right) &\geq& 0.
\end{eqnarray}
\end{subequations}
Addition these equations in (\ref{sm}) gives us
\begin{equation}
\left( J^{\;'}(u_h) - J^{\;'}(u), u-u_h \right) + \big( \zeta, u-u_h \big) \geq 0.
\end{equation}
By the mean value theorem, we obtain
\begin{equation}
-J^{\;''}(\widehat{u}) (u-u_h)^2 + (\zeta, u-u_h) \geq 0
\end{equation}
for some $\widehat u \in \{ v \in U_{ad}:\; v=u+t(u_h -u), \;\; t\in(0,1)\}$. Then, when we apply the SSC (\ref{ssc}) with Cauchy-Schwarz inequality, we find
\begin{equation}
\delta^{\;'} \|u_h -u\|^{2}_{L^2(Q_T)} \leq \| \zeta\|_{L^2(Q_T)} \|u_h -u\|_{L^2(Q_T)},
\end{equation}
which is the desired result.
\end{proof}
Here we follow  \cite[Remark~3.3]{EKammann_FTroltzsch_SVolkwein} and select $\delta^{\;'}:=\delta /2$. If $u_h$  belongs to the neighborhood  of $u$, we can estimate in the following
\begin{equation}\label{sscest}
\|u - u_h\|_{L^2(Q_T)} \leq \frac{2}{\delta^{\;'}} \|\zeta\|_{L^2(Q_T)}.
\end{equation}

\begin{remark}
We also remark  that the presence of the so-called Tikhonov parameter $\frac{\omega_u}{2}$ in the cost functional (\ref{M4}) is extremely important. The standard second-order optimality conditions do not hold for vanishing Tikhonov parameter $\omega_u=0$. By introducing new critical cones, the second-order optimality conditions are established in \cite{ECasas_CRyll_FTroltzsch_2015} for sparse optimal control governed by FitzHugh-Nagumo system, and in  \cite{ECasas_FTroltzsch_2016} for general nonlinear functions. Hence, the new second order conditions are used  for proving the stability of locally optimal solutions with respect to  $\omega_u \to 0$. This theoretical result was confirmed for one dimensional wave solutions of the classical FHN equation with sparse controls in \cite{ECasas_CRyll_FTroltzsch_2015,CRyll_JLober_SMartens_HEngel_FTroltzsch_2016}. In this paper, we apply the theory introduced in \cite{ECasas_CRyll_FTroltzsch_2015,ECasas_FTroltzsch_2016}
for two dimensional wave solutions of the convective FHN equations.
\end{remark}

%%%%%%%%%%%%%%%%%%%%%%%%%%%%%%%%%%%%%%%%%%%%%%%%%%%%%%%%%%%%%%%%%%%%%%%%%%%%%%%%%%%%%%%%%%%%
\section{Discretization in space and time}
\label{disc_state}

We discretize the  state system  (\ref{M1}) and the adjoint system  (\ref{M9}) by  the symmetric interior penalty Galerkin (SIPG) in space and  the backward  Euler discretization in time.

\subsection{SIPG discretization of the state and adjoint equations in space}

%For the space discretization we choose symmetric interior penalty Galerkin (SIPG) method for the diffusion and an upwind discretization for the convection  due to the its symmetric property.
 The interior penalty dG methods are well established in the literature and the details can be found in the classical texts like
 \cite{DNArnold_FBrezzi_BCockburn_LDMarini_2002a,BRiviere_2008a}.

We assume that the domain $\Omega$ is polygonal domain. We denote $\{ \mathcal{T}_h\}_h$ as a family of shape-regular simplicial triangulations of $\Omega$, see, e.g., \cite{PGCiarlet_2002}. Each mesh $\mathcal{T}_h$ consists of closed triangles such that $\overline{\Omega} = \bigcup_{K \in \mathcal{T}_h} \overline{K}$ holds. We assume that the mesh is regular in the following sense: for different triangles
$K_i, K_j \in \mathcal{T}_h$, $i \not= j$, the intersection  $K_i \cap K_j$ is either empty or a vertex or an edge, i.e., hanging nodes are not allowed. The diameter of an element $K$ and the length of an edge $E$ are denoted by $h_{K}$ and $h_E$, respectively.

We split the  set of all edges $\mathcal{E}_h$ into the set $\mathcal{E}^0_h$ of interior edges, the set $\mathcal{E}^{D}_h$ of Dirichlet boundary edges and  the set $\mathcal{E}^{N}_h$  of Neumann boundary edges so that $\mathcal{E}_h=\mathcal{E}^{B}_h  \cup \mathcal{E}^{0}_h$ with $\mathcal{E}^B_h=\mathcal{E}^{D}_h \cup \mathcal{E}^{N}_h$.
Let $\mathbf{n}$ denote the unit outward normal to $\partial \Omega$. For the activator $y$ and the inhibitor $z$, we define the
inflow and outflow parts of $\partial \Omega$ by
\[
          \Gamma^- = \set{x \in \partial \Omega}{ \mathbf{V}(x) \cdot \mathbf{n}(x) < 0}, \quad
          \Gamma^+ = \set{x \in \partial \Omega}{ \mathbf{V}(x) \cdot \mathbf{n}(x) \geq 0}.
\]
Similarly, the inflow and outflow boundaries of an element $K$ are defined by
\[
\partial K^-=\set{x \in \partial K}{\mathbf{V}(x) \cdot \mathbf{n}_{K}(x) <0}, \quad \partial K^+=\set{x \in \partial K}{\mathbf{V}(x) \cdot \mathbf{n}_{K}(x) \geq 0},
\]
where $\mathbf{n}_{K}$ is the unit normal vector on the boundary $\partial K$ of an element $K$.

Let the edge $E$ be a common edge for two elements $K$ and $K^e$. For a piecewise continuous scalar function $y$, there are two traces of $y$ along $E$, denoted by $y|_E$ from inside $K$ and $y^e|_E$ from inside $K^e$. Then, the jump and average of $y$ across the edge $E$ are defined by:
\begin{align}
\jump{y}&=y|_E\mathbf{n}_{K}+y^e|_E\mathbf{n}_{K^e}, \quad
\average{y}=\frac{1}{2}\big( y|_E+y^e|_E \big).
\end{align}
Similarly, for a piecewise continuous vector field $\nabla y$, the jump and average across an edge $E$ are given by
\begin{align}
\jump{\nabla y}&=\nabla y|_E \cdot \mathbf{n}_{K}+\nabla y^e|_E \cdot \mathbf{n}_{K^e}, \quad
\average{\nabla y}=\frac{1}{2}\big(\nabla y|_E+\nabla y^e|_E \big).
\end{align}
For a boundary edge $E \in K \cap \partial\Omega$, we set $\average{ y}= y$ and $\jump{y}=y\mathbf{n}$.

Recall that for dG  methods, we do not impose continuity constraints on the trial and test functions across the element interfaces. As a consequence, the weak formulation must include jump terms across interfaces, and typically penalty terms are added to control the jump terms. Then, we define the discontinuous discrete space as follows:
\begin{align}\label{DG1}
W_h = \set{w \in L^2(\Omega)}{w \mid_{K}\in \mathbb{P}^1(K), \quad \forall K \in \mathcal{T}_h},
\end{align}
where $\mathbb{P}^1(K)$  is the set of piecewise linear  polynomials defined on $K$. We note that the space of discrete states and the space of test functions are identical due to the weak treatment of boundary conditions for dG methods.

Then, the semi-discrete formulation of the state system (\ref{M1})  for $\forall w \in W_h$  and $ t \in (0,T]$ becomes
\begin{subequations}\label{DG2}
\begin{eqnarray}
\hspace{-5mm}\left ( \frac{d y_h}{dt}, w \right) + a_{h,y}(y_h,w) + b_{h,y}(y_h,w) + c_{h,z}(z_h,w) &=& \ell_{h,y}(w) + (u_h,w), \\
\hspace{-5mm} (y_h(\cdot,0), w) &=& (y_0,w), \\
\hspace{-5mm}\left ( \frac{d z_h}{dt}, w \right) + a_{h,z}(z_h,w) + b_{h,z}(z_h,w) + c_{h,y}(y_h,w) &=& \ell_{h,z}(w), \\
\hspace{-5mm} (z_h(\cdot,0), w) &=& (z_0,w),
\end{eqnarray}
\end{subequations}
where the (bi)-linear terms are defined for $i=y,z$ and $\forall w \in W_h$
\begin{subequations}\label{DG3}
\begin{align}
a_{h,i}(v,w)=& \sum \limits_{K \in \mathcal{T}_h} \int \limits_{K} d_i \nabla v \cdot  \nabla w \; dx \nonumber \\
           & -\sum \limits_{E \in \mathcal{E}^0_h \cup \mathcal{E}^D_h} \int \limits_E \Big( \average{d_i \nabla v} \cdot \jump{w} +
               \average{d_i \nabla w} \cdot \jump{v} \Big)  \; ds  \nonumber \\
           &+ \sum \limits_{E \in \mathcal{E}^0_h \cup \mathcal{E}^D_h} \frac{\sigma d_i}{h_E} \int \limits_E \jump{v} \cdot \jump{w} \; ds +  \sum \limits_{K \in \mathcal{T}_h} \int \limits_{K} \mathbf{V} \cdot \nabla v w  \; dx \nonumber \\
           &+ \sum \limits_{K \in \mathcal{T}_h}\; \int \limits_{\partial K^{-} \backslash \partial \Omega} \mathbf{V} \cdot \mathbf{n} (v^e-v)w \; ds - \sum \limits_{K \in \mathcal{T}_h} \; \int \limits_{\partial K^{-} \cap \Gamma^{-}} \mathbf{V} \cdot \mathbf{n} v w  \; ds, \\
\ell_{h,i}(w)=&  \sum \limits_{E \in \mathcal{E}^D_h} \int \limits_E  i_{D} \Big(\frac{\sigma d_i}{h_E} \mathbf{n} \cdot \jump{w}
              - \average{d_i \nabla w}  \Big)\; ds \nonumber \\
           & - \sum \limits_{K \in \mathcal{T}_h} \; \int \limits_{\partial K^{-} \cap \Gamma^{-}} \mathbf{V} \cdot \mathbf{n} \; i_{D} w  \; ds,
\end{align}
\begin{align}
b_{h,y}(y,w) =&\sum \limits_{K \in \mathcal{T}_h} \int \limits_{K} g(y)  w  \; dx,  \qquad
c_{h,y}(y,w) =\sum \limits_{K \in \mathcal{T}_h} \int \limits_{K} -\epsilon c_3 y w \; dx,  \\
b_{h,z}(z,w) =&\sum \limits_{K \in \mathcal{T}_h} \int \limits_{K} \epsilon z w \; dx,  \qquad \quad
c_{h,z}(z,w) =\sum \limits_{K \in \mathcal{T}_h} \int \limits_{K} z  w \; dx,
\end{align}
\end{subequations}
where the parameter $\sigma \in \mathbb{R}_0^+$ is called the penalty parameter which should be sufficiently large to ensure the stability of the dG discretization; independent of the mesh size $h$ and of the diffusion coefficients $d_i$ as described in \cite[Sec.~2.7.1]{BRiviere_2008a} with a lower bound depending only on the polynomial degree. Large penalty parameters  decrease the jumps across element interfaces, which can affect the numerical approximation. However, the dG approximation can converge to the continuous Galerkin approximation as the penalty parameter goes to infinity (see, e.g., \cite{ACangiani_JChapman_HGeorgoulis_MJensen_2014} for details).

For each time step, we can expand the discrete solutions of the states $y, z$, and the control $u$ as
\begin{align}\label{DG4}
y_h(t) = \sum \limits_{i=1}^{N} \sum \limits_{j=1}^{n_{k}} \vec{y}_{j}^{\,i} \phi_{j}^{\,i}, \quad
z_h(t) = \sum \limits_{i=1}^{N} \sum \limits_{j=1}^{n_{k}} \vec{z}_{j}^{\,i} \phi_{j}^{\,i}, \quad \hbox{and} \quad
u_h(t) = \sum \limits_{i=1}^{N} \sum \limits_{j=1}^{n_{k}} \vec{u}_{j}^{\,i} \phi_{j}^{\,i},
\end{align}
where $\vec{y}_{j}^{\,i}$, $\vec{z}_{j}^{\,i}$, $\vec{u}_{j}^{\,i}$,  and $\phi_{j}^{\,i}$ are the unknown coefficients and the basis functions, respectively, for $j=1,2,\cdots, n_{k}$ and $i=1,2, \cdots, N$. The number $N$ denotes the number of dG elements and $n_{k}=(k+1)(k+2)/2$ is the local dimension of each dG element depending on the order $k$ of the polynomial basis.

Inserting (\ref{DG4}) into (\ref{DG2}), we obtain
\begin{subequations}\label{DG5}
\begin{eqnarray}
 \mathbf{M} \frac{d \vec{y}}{dt} +  \mathbf{A}_y \vec{y}  + \mathbf{b}_y(\vec{y}) + \mathbf{C}_z \vec{z}  &=&  \mathbf{l}_y + \mathbf{M} \vec{u}, \\
 \mathbf{M} \frac{d \vec{z}}{dt} +  \mathbf{A}_z \vec{z}  + \mathbf{B}_z \vec{z} + \mathbf{C}_y \vec{y} &=&  \mathbf{l}_z,
\end{eqnarray}
\end{subequations}
where $\vec{y}$, $\vec{z}$ and $\vec{u}$ are the unknown coefficient vectors, i.e.,
\begin{eqnarray*}
\vec{y}&=&(\vec{y}_{1}^{\,1}, \cdots, \vec{y}^{\,1}_{n_k}, \cdots, \vec{y}^{\,N}_{1}, \cdots, \vec{y}^{\,N}_{n_k}), \\
\vec{z}&=&(\vec{z}_{1}^{\,1}, \cdots, \vec{z}^{\,1}_{n_k}, \cdots, \vec{z}^{\,N}_{1}, \cdots, \vec{z}^{\,N}_{n_k}), \\
\vec{u}&=&(\vec{u}_{1}^{\,1}, \cdots, \vec{u}^{\,1}_{n_k}, \cdots, \vec{u}^{\,N}_{1}, \cdots, \vec{u}^{\,N}_{n_k}),
\end{eqnarray*}
$\mathbf{M}$ is the mass matrix, $\mathbf{A}_y$ and $\mathbf{A}_z$ are the stiffness matrices corresponding to $a_h(y_h,w)$ and $a_h(z_h,w)$, respectively, and $\mathbf{l}_y$ and $\mathbf{l}_z$ are vectors corresponding to $\ell_{h,y}(w)$ and $\ell_{h,z}(w)$, respectively. Moreover, $\mathbf{C}_y$, $\mathbf{B}_z$, and $\mathbf{C}_z$ are the matrices corresponding to $c_{h,y}$, $b_{h,z}$, and $c_{h,z}$, respectively. Further, $\mathbf{b}_y(\vec{y})$ is the vector corresponding to the nonlinear part  $b_{h,y}(y,w)$.

In the literature, there exist two approaches for the solution the OCP (\ref{M4})
numerically, i.e., \emph{discretize-then-optimize} (DO)  and \emph{optimize-then-discretize} (OD). In the DO approach, the objective function (\ref{M4}) and the state equation (\ref{M1}) are discretized first, and then the discrete optimality
system is formed. In the OD approach, the optimality conditions consisting of
the state system (\ref{M1}), the adjoint system (\ref{M9}) and the variational form (\ref{M8}) are derived. Then, the infinite dimensional optimality system is discretized. We here apply the \emph{optimize-then-discretize} approach. The discretization and optimization commute for the SIPG discretized linear steady state convection-diffusion-reaction OCPs \cite{HYucel_MHeinkenschloss_BKarasozen_2013,HYucel_BKarasozen_2014}. For linear time dependent OPCs the OD and DO approaches commute for the dG discretization in time and the SIPG discretization in space \cite{Akman14}. The commutativity property is no more valid for semi-linear elliptic and parabolic OCPs.

In the \emph{optimize-then-discretize} approach, we derive directly semi-discrete form of the adjoint equations (\ref{M9}) using the SIPG discretization as for the state equations (\ref{DG2}):
\begin{subequations}\label{N0}
\begin{align*}
\left ( -\frac{d p_h}{dt}, w \right) + a_{h,p}(p_h,w) + b_{h,p}(p_h,w) + c_{h,q}(q_h,w) =& \omega_{Q}^y
\big(y_h,w \big) + \ell_{h,p}(w), \\
 (p_h(\cdot,T),w) =&  \omega_T^y \left( \big( y_h(\cdot, T)-y_T(x) \big),w \right), \\
\left ( -\frac{d q_h}{dt}, w \right) + a_{h,q}(q_h,w) + b_{h,q}(q_h,w) + c_{h,p}(p_h,w) =& \omega_{Q}^z \big(z_h,w \big) + \ell_{h,q}(w), \\
 (q_h(\cdot,T),w) =&  \omega_T^z \left( \big( z_h(\cdot, T)-z_T(x) \big),w \right),
\end{align*}
\end{subequations}
where the bilinear forms $a_{h,p}$ and $a_{h,q}$ are similar to the state ones (\ref{DG3}) with the negative convection terms, i.e., $-\mathbf{V}(x,y)$. As an extra term, they just contain the contribution of the mixed boundary conditions, i.e.,
\[
\sum \limits_{E \in \mathcal{E}^N_h}  \int \limits_E \big( \mathbf{V} \cdot \mathbf{n} \big) p_h w\; ds \quad \hbox{and} \quad \sum \limits_{E \in \mathcal{E}^N_h}  \int \limits_E \big( \mathbf{V} \cdot \mathbf{n} \big) q_h w\; ds.
\] The other terms are given by
\begin{subequations}\label{N00}
\begin{align}
b_{h,p}(p,w) &= \sum \limits_{K \in \mathcal{T}_h} \int \limits_{K} g_y(y) p w  \; dx, \qquad
c_{h,q}(q,w) = \sum \limits_{K \in \mathcal{T}_h} \int \limits_{K} -\epsilon c_3 q w \; dx,
\end{align}
\begin{align}
b_{h,q}(q,w) &= \sum \limits_{K \in \mathcal{T}_h} \int \limits_{K} \epsilon q w  \; dx, \qquad \qquad
c_{h,p}(p,w) = \sum \limits_{K \in \mathcal{T}_h} \int \limits_{K} p w \; dx,  \\
\ell_{h,p}(w) &= -\omega_Q^y \sum \limits_{K \in \mathcal{T}_h} \int \limits_{K} y_Q w \; dx, \qquad
\ell_{h,q}(w) = -\omega_Q^z \sum \limits_{K \in \mathcal{T}_h} \int \limits_{K} z_Q w \; dx.
\end{align}
\end{subequations}

%Analogously to $y$, $z$ and $u$, we discretize $p$ and $q$ by
%\begin{align}\label{N1}
%p_h(t) = \sum \limits_{i=1}^{N} \sum \limits_{j=1}^{n_{k}} \vec{p}_{j}^{\,i} \phi_{j}^{\,i} \quad \hbox{and} \quad
%q_h(t) = \sum \limits_{i=1}^{N} \sum \limits_{j=1}^{n_{k}} \vec{q}_{j}^{\,i} \phi_{j}^{\,i}.
%\end{align}

%%%%%%%%%%%%%%%%%%%%%%%%%%%%%%%%%%%%%%%%%%%%%%%%%%%%%%%%%%%%%%%%%%%%%%%%%%%%

\subsection{Time Discretization}
\label{disc_opg}

%After  spatial discretization of the state system (\ref{M1}), we now ready for temporal discretization.

We solve the nonlinear system of ordinary differential equations (\ref{DG5}) in time by  the backward  Euler method. We  first divide the time interval $[0,T]$ $
0=t_0 < t_1 < \ldots < t_{N_T}=T$,
with  step size $\tau_n = t_{n}-t_{n-1}, \;\; n=1,2,\ldots, N_T$ and then obtain the following fully discretized state system:
\begin{subequations}\label{DG6}
\begin{eqnarray}
 \frac{1} {\tau_n} \mathbf{M} \big( \vec{y}_n - \vec{y}_{n-1} \big)  +  \mathbf{A}_y \vec{y}_{n}  + \mathbf{b}_y(\vec{y}_{n}) + \mathbf{C}_z \vec{z}_{n}  &=&  \mathbf{l}_y^n + \mathbf{M} \vec{u}_n, \\
 \frac{1} {\tau_n} \mathbf{M} \big( \vec{z}_n - \vec{z}_{n-1} \big) +  \mathbf{A}_z \vec{z}_n  + \mathbf{B}_z \vec{z}_n + \mathbf{C}_z \vec{y}_{n} &=&  \mathbf{l}_z^n
\end{eqnarray}
\end{subequations}
for $n=1,2,\ldots, N_T$. We solve the nonlinear system  by Newton's method.

%Then, by proceeding as the state system done in section~\ref{disc_state}, we obtain the following
We discretize the adjoint variables $p$ and $q$ analogously to the state variables $y$, and $z$ as
\begin{align}\label{N1}
p_h(t) = \sum \limits_{i=1}^{N} \sum \limits_{j=1}^{n_{k}} \vec{p}_{j}^{\,i} \phi_{j}^{\,i} \quad \hbox{and} \quad
q_h(t) = \sum \limits_{i=1}^{N} \sum \limits_{j=1}^{n_{k}} \vec{q}_{j}^{\,i} \phi_{j}^{\,i}.
\end{align}
By inserting (\ref{N1}) into (\ref{N0}) and applying the backward Euler discretization in time, we obtain the following  fully discretized adjoint system:
\begin{subequations}\label{N2}
\begin{align}
 \frac{1} {\tau_n} \mathbf{M} \big( \vec{p}_{n-1} - \vec{p}_{n} \big)  +  \mathbf{A}_p \vec{p}_{n-1}  + \mathbf{B}_p(\vec{y}_{n}) \vec{p}_{n-1} + \mathbf{C}_q \vec{q}_{n-1}  =&  \omega_Q^y \big( \mathbf{M} \vec{y}_{n}+ \mathbf{l}_p^n \big), \\
 \frac{1} {\tau_n} \mathbf{M} \big( \vec{q}_{n-1} - \vec{q}_{n} \big) +  \mathbf{A}_q \vec{q}_{n-1}  + \mathbf{B}_q \vec{q}_{n-1} + \mathbf{C}_p \vec{p}_{n-1} =& \omega_Q^z \big( \mathbf{M} \vec{z}_{n}+ \mathbf{l}_q^n \big)
\end{align}
\end{subequations}
for $n= N_T, \ldots,2,1$. The  matrices  $\mathbf{A}_p, \mathbf{A}_q$ correspond to the discretized  bilinear forms $a_{h,p}$ and $a_{h,q}$, respectively, whereas the matrices $\mathbf{B}_p(y)$, $\mathbf{C}_p$, $\mathbf{B}_q$, $\mathbf{C}_q$ and the vectors $\mathbf{l}_p, \mathbf{l}_q$  are the discretized forms in (\ref{N00}), respectively.

%We can now solve the optimal control problem governed by the convective FitzHugh-Nagumo system (\ref{M4}).

%%%%%%%%%%%%%%%%%%%%%%%%%%%%%%%%%%%%%%%%%%%%%%%%%%%%%%%%%%%%%%%%%%%%%%%%%
\section{Numerical results}
\label{numeric}

In this section, we provide numerical results for the OCP governed by the convective FHN equation (\ref{M4}).
%As a optimization procedure, we prefer to apply the nonlinear conjugate approach (namely, the one Polak-Ribiere) along with the approximate Wolfed step-size rule for the linear search to guarantee the descent direction.
There exists several optimization algorithms for OCPs governed by semi-linear equations. We have used in this paper the  projected  nonlinear conjugate gradient (CG) method  \cite{WWHager_HZhang_2006aa,WWHager_HZhang_2006}, which was applied to the Schl\"ogl and FHN equations \cite{RBuchholz_HEngel_EKammann_FTroltzsch_2013,ECasas_CRyll_FTroltzsch_2013,ECasas_CRyll_FTroltzsch_2015}. The projected nonlinear CG algorithm is outlined below:
\begin{description}
  \item[Initialization] Choose an initial guess $u_0$, an initial step size $s_0$ and  stopping tolerances $\hbox{Tol}_1$ and $\hbox{Tol}_2$. Then, compute
         \begin{itemize}
           \item initial states $(y_0,z_0)=(y_{u_{0}},z_{u_{0}})$ by solving (\ref{DG6}),
           \item initial adjoints $(p_0,q_0)=(p_{y_{0},z_{0}}, q_{y_{0},z_{0}})$ by solving (\ref{N2}),
           \item subgradient of $j$ by (\ref{M11}),  i.e., $\lambda_0 = \lambda_{u_{0}, p_{0}}$,
           \item subgradient of $J$,  i.e, $g_0=\omega_c u_0 + p_0 + \mu \lambda_0$,
           \item anti-subgradient of $J$, i.e., $d_0=-g_0$.
         \end{itemize}
         Set k:=0.
  \item[Step 1. (New subgradients)] Update
         \begin{itemize}
           \item control, i.e., $u_{k+1}=u_k + s_k d_k$,
           \item states  $(y_{k+1},z_{k+1})=(y_{u_{k+1}},z_{u_{k+1}})$ by solving (\ref{DG6}),
           \item adjoints $(p_{k+1},q_{k+1})=(p_{y_{k+1},z_{k+1}}, q_{y_{k+1},z_{k+1}})$  by solving (\ref{N2}),
           \item subgradient of $j$ by (\ref{M11}), i.e., $\lambda_{k+1} = \lambda_{u_{k+1}, p_{k+1}}$,
           \item subgradient of $J$, i.e, $g_{k+1}=\omega_c u_{k+1} + p_{k+1} + \mu \lambda_{k+1}$.
         \end{itemize}
  \item[Step 2. (Stopping Criteria)] Stop if $\|g_{k+1}\| < \hbox{Tol}_1$ or $\|J_{k+1}-J_k\| \leq \hbox{Tol}_2$.
  \item[Step 3. (Direction of descent)] Compute the conjugate direction $\beta_{k+1}$ according to one of the update formulas such as Hestenes-Stiefel, Polak-Ribiere, Fletcher-Reeves, and Hager-Zhang, see e.g., \cite{RBuchholz_HEngel_EKammann_FTroltzsch_2013,WWHager_HZhang_2006} for details.
      \[
         d_{k+1} = -g_{k+1} + \beta_{k+1} d_k.
      \]
   \item[Step 4. (Step size)] Select step size $s_{k+1}$ according to some standard options such as bisection, strong Wolfed-Powell, see e.g., \cite{RBuchholz_HEngel_EKammann_FTroltzsch_2013,WWHager_HZhang_2006aa} for details. Set $k=:k+1$ and go to Step 1.
\end{description}

For the computation of the reduced Hessian, we use the BFGS algorithm \cite{JNocedal_SJWright_2006a,RHerzog_KKunisch_2010}:
\begin{itemize}
  \item Set $H_{0}=I$.
  \item Update for $k=1,2,\ldots$
  \[
    H_{k+1} = H_{k} + \frac{q_{k} q_{k}^{T}}{q_{k}^{T} r_{k}} - \frac{\big(H_{k} r_{k}\big) \big(H_{k} r_{k}\big)^{T}} {r_{k}^{T} H_{k} s_{k}},
  \]
  where  $r_{k}= u_{k+1} -u_{k}$ and $q_{k} = g_{k+1} -g_{k}$.
\end{itemize}
After setting the reduced Hessian matrix at each time level, we compute the smallest eigenvalue of the Hessian by using MATLAB\textsuperscript{\textregistered} \texttt{eig}.

In all numerical experiments, we started with the initial control $u=0$. The optimization algorithm is terminated, when the gradient is less than $\hbox{Tol}_1=10^{-3}$ or the difference between the successive cost functionals is less than $\hbox{Tol}_2=10^{-5}$. We use uniform step sizes in time  $\Delta t=0.05$. The penalty parameter $\sigma$ in (\ref{DG3}) is chosen as  $\sigma=6$ on the interior edges and $\sigma=12$  on the boundary edges.

If it is not specified, we use the following parameters in all numerical examples
\[
c_1=9, \quad c_2 = 0.02, \quad c_3 = 5, \quad \epsilon = 0.1, \quad d_{1}=d_{2}=1
\]
on  a rectangular box $\Omega=[0,L] \times [0,H]$ with $L=100$ and $H=5$ as in  \cite{EAErmakova_EEShnol_MAPanteleev_AAButylin_VVolpert_FIAtaullakhanov_2009}. Further, the final time is taken as $T=1$.

\indent All simulations are generated on Intel(R) Core(TM)  i7-4720HQ CPU@2.60GHz, 16GB RAM, Windows 8, with MATLAB\textsuperscript{\textregistered} R2014a (64-bit).

%%%%%%%%%%%%%%%%%%%%%%%%%%%%%%%%%%%%%%%%%%%%%%%%%%%%%%%%%%%%%%%%%

\subsection{Optimal control in the space-time domain}\label{Ex1}
We first consider  the OCP with the desired state functions defined  in the whole space-time domain $Q_T$ with $\omega_{Q}^{y}=\omega_{Q}^{z}=1$, $\omega_{T}^{y}=\omega_{T}^{z}=0$,  and $\omega_u=10^{-5}$. The step sizes in space are taken as $\Delta x_{1}=\Delta x_{2} =0.5$. The desired states are chosen as the solution of uncontrolled FHN equation (\ref{M1})
\[
y_Q(x,t) = \left\{
             \begin{array}{ll}
               y_{\mathrm{ nat}}(x,t), & \hbox{if} \;\; t \leq T/2,  \\
               0, & \hbox{otherwise},
             \end{array}
           \right.
\qquad
z_Q(x,t) = \left\{
             \begin{array}{ll}
               z_{\mathrm{nat}}(x,t), & \hbox{if} \;\; t \leq T/2,  \\
               0, & \hbox{otherwise},
             \end{array}
           \right.
\]
with the initial conditions given by
\[
y_0(x,t) = \left\{
             \begin{array}{ll}
               0.1, & \hbox{if} \;\; 0 \leq x_1 \leq 0.1,  \\
               0, & \hbox{otherwise},
             \end{array}
           \right.
\qquad
z_0(x,0) = 0.
\]
Further, we define the admissible set of controls as
\[
\mathcal{U}_{ad} := \{ u \in L^{\infty}(Q): \; -0.2 \leq u(x,t) \leq 0 \;\; \hbox{for   a.e   } (x,t) \in Q \}.
\]

\begin{figure}[H]
   \includegraphics[width=1\textwidth]{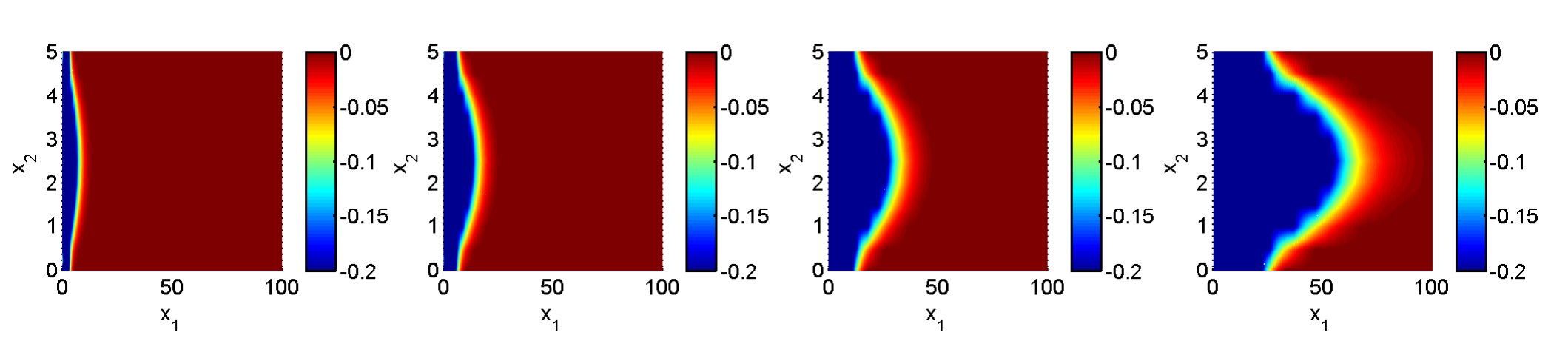}
   \includegraphics[width=1\textwidth]{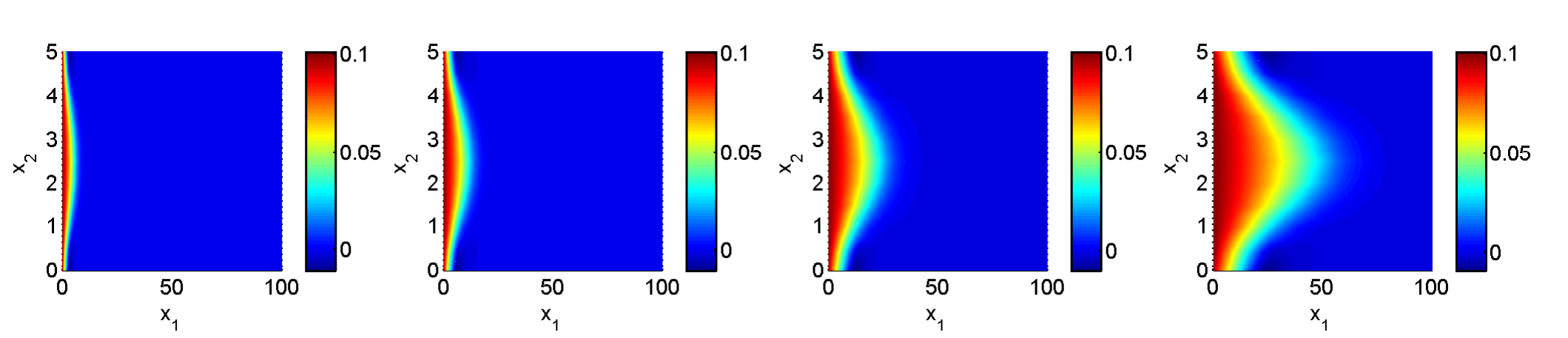}
    \caption{Example~\ref{Ex1}: optimal controls $u$ (top) and associated states $y$  (bottom) at $t=0.75$  without sparse control for  $V_{\max}$=16, 32, 64, 128 (from left to right).}
    \label{Fig:Ex1_const}
\end{figure}

\begin{table}[H]
\centering
\begin{tabular}{r|c|r|c|r}
  $V_{max}$ & $J$ &  $\#$ite.  & $\#$search  & $\#$Newton    \\
\hline
  $16$  & 2.91e-2 & 64     & 245  &  787   \\
  $32$  & 5.80e-2 & 77     & 297  &  956   \\
  $64$  & 1.16e-1 & 96     & 373  &  1203   \\
  $128$  & 2.30e-1 & 123     & 481  & 1554   \\
\hline
\end{tabular}
\caption{Example~\ref{Ex1}: cost functional $J$, number of nonlinear CG iterations, line searches, Newton steps without sparse control.}
\label{Ex1_table1}
\end{table}

\begin{figure}[H]
\centering
    \subfloat[Control $u$]{%
      \includegraphics[width=0.35\textwidth]{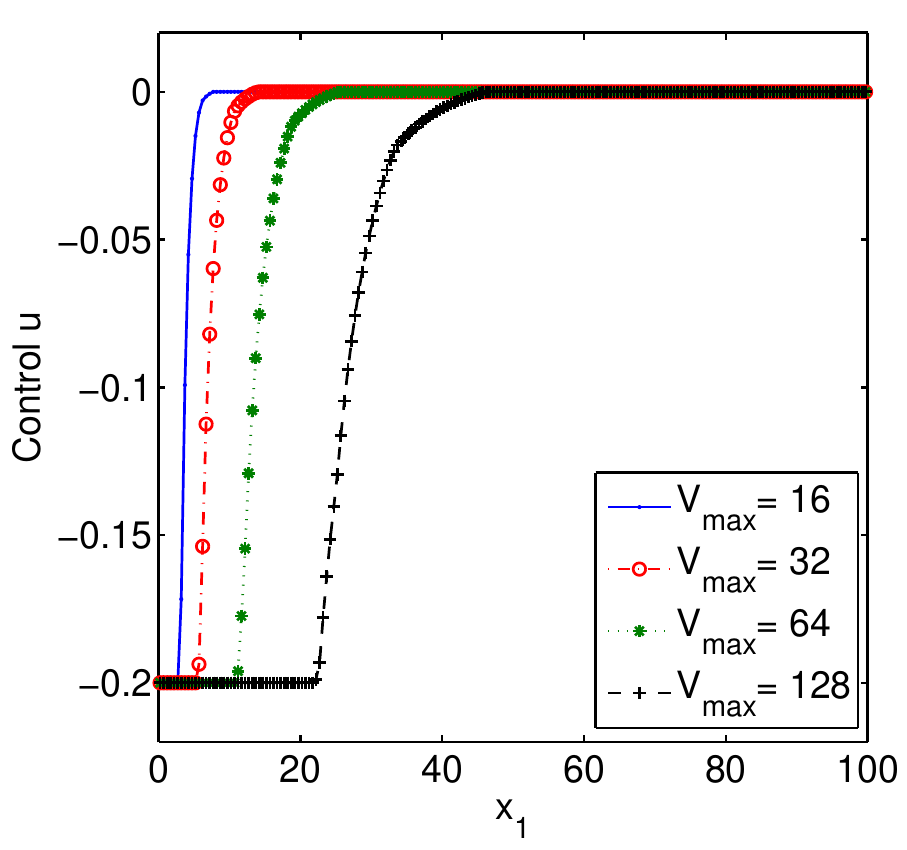}
    }
    \qquad \quad
    \subfloat[State $y$]{%
      \includegraphics[width=0.35\textwidth]{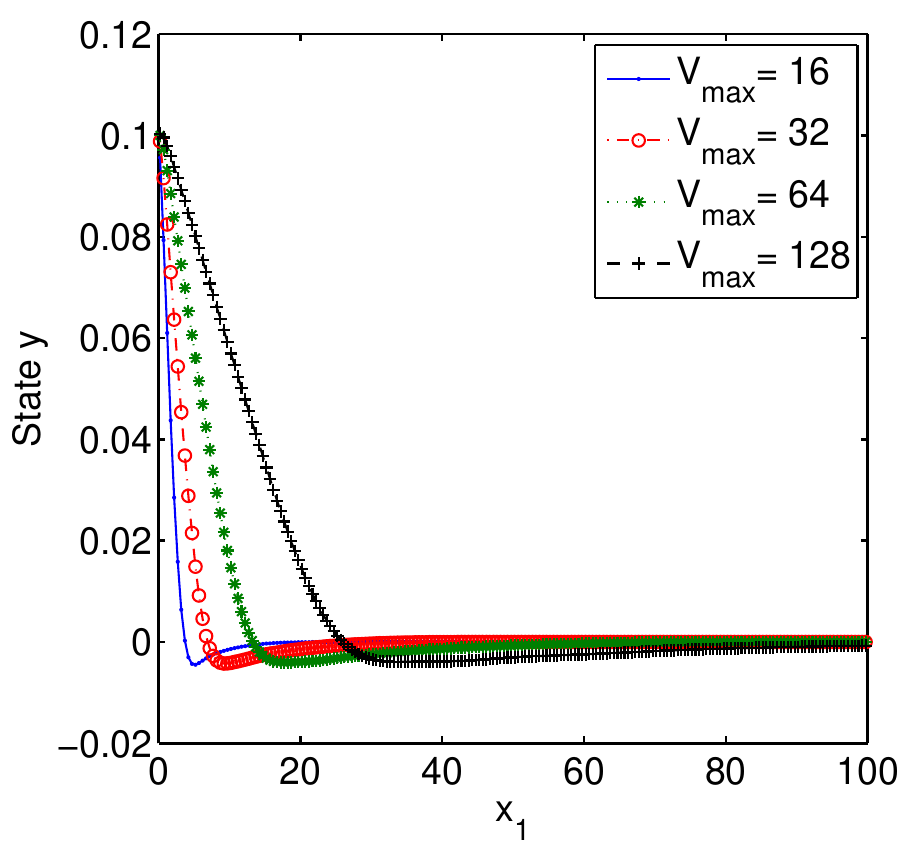}
    }
    \caption{Example~\ref{Ex1}: profiles of optimal controls $u$ and associated states $y$ along $x_1$ axis at $t=0.75$ without sparse control for various values of $V_{max}$.}
    \label{Fig:Ex1_const_Vmax}
\end{figure}

\begin{figure}[H]
\centering
    \subfloat[Control $u$]{%
      \includegraphics[width=0.35\textwidth]{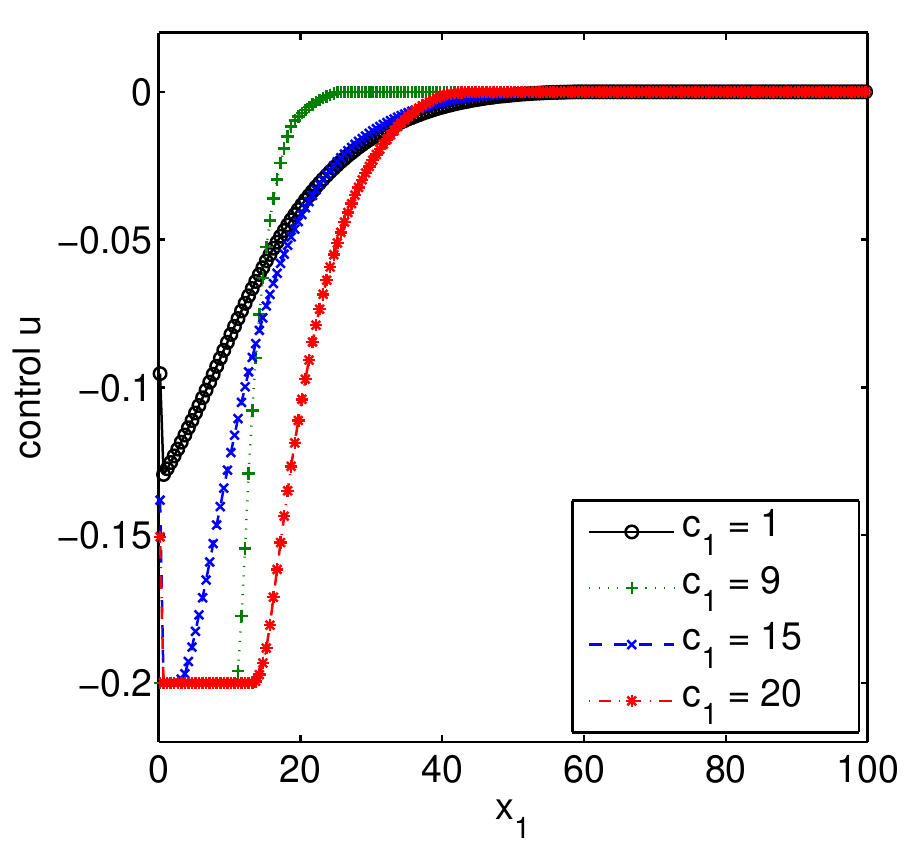}
    }
    \qquad \quad
    \subfloat[State $y$]{%
      \includegraphics[width=0.35\textwidth]{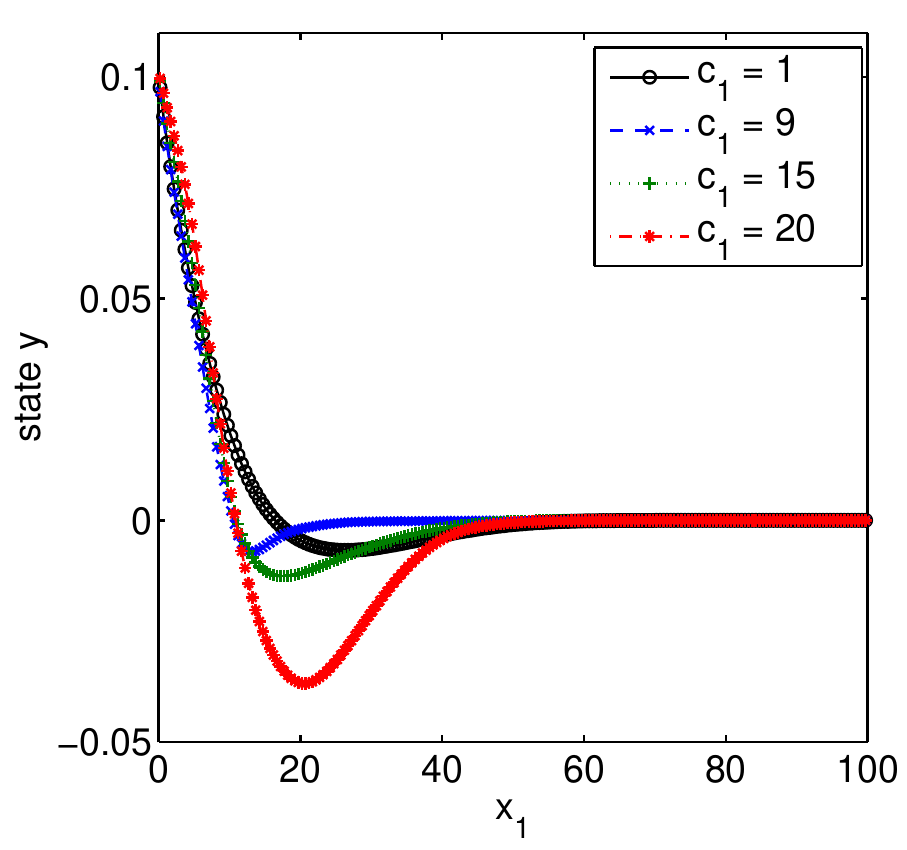}
    }
    \caption{Example~\ref{Ex1}: profiles of optimal controls $u$ and associated states $y$ along $x_1$ axis for  $V_{\max}=64$  without sparse control for various values of $c_1$.}
    \label{Fig:Ex1_const_c1}
\end{figure}

\begin{figure}[H]
   \includegraphics[width=1\textwidth]{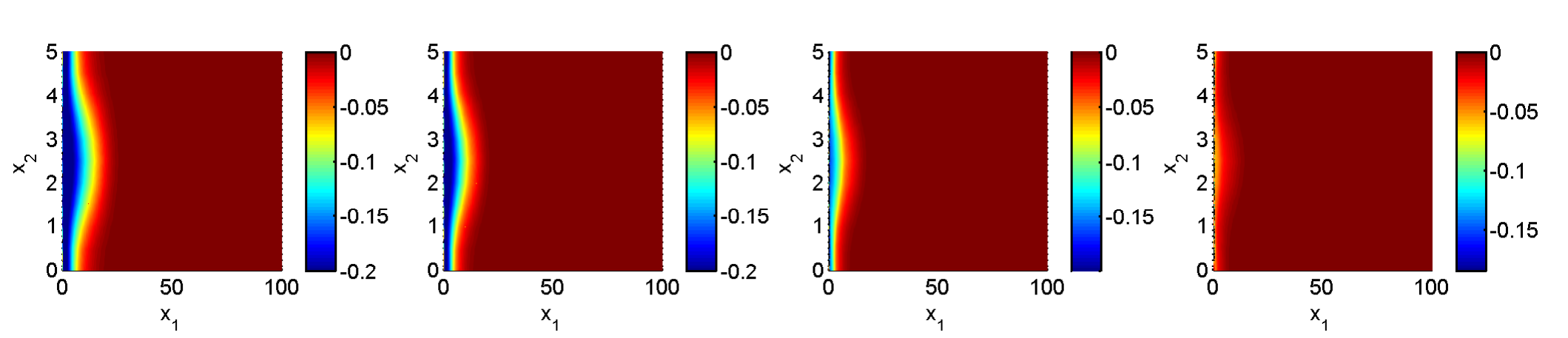}
   \includegraphics[width=1\textwidth]{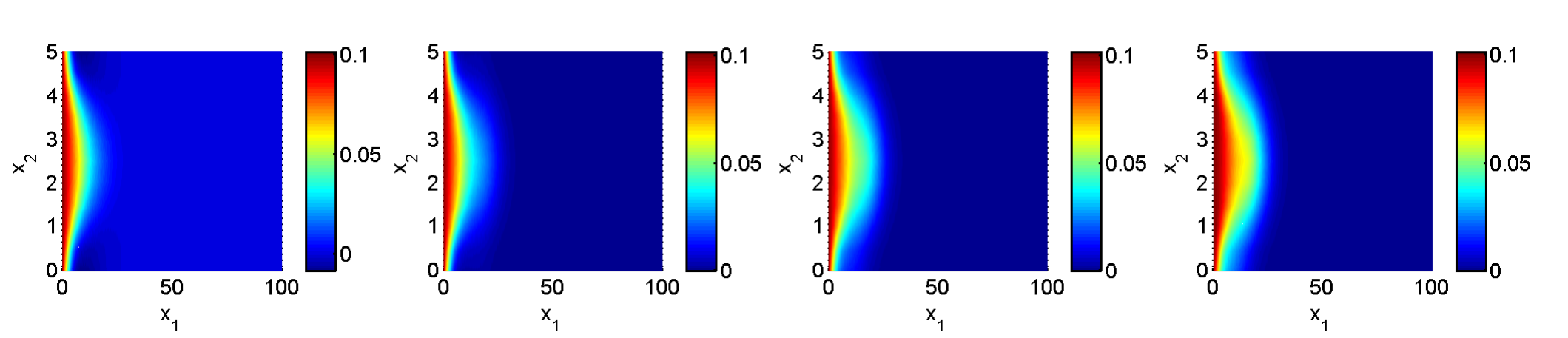}
    \caption{Example~\ref{Ex1}: optimal controls $u$ (top) and associated states $y$ (bottom) for the sparse parameters $\mu=1/500, 1/100, 1/50$,  $1/35$ (from left to right) and  $V_{\max}=32$ at $t=0.75$.}
    \label{Fig:Ex1_spr}
\end{figure}

We first investigate the numerical solutions of the optimization problem (\ref{M4}) without the sparse parameter, i.e.,  $\mu=0$. Figure~\ref{Fig:Ex1_const} demonstrates the computed solutions of the control $u$ and their associated  states $y$ at $t=0.75$ for various values of the $V_{\max}$, respectively. The control $u$, bounded by the box constraints, exhibits the same behavior of the state $y$. When the value of $V_{\max}$ increases, the controlled solutions become more curved, as for the uncontrolled solutions \cite{EAErmakova_EEShnol_MAPanteleev_AAButylin_VVolpert_FIAtaullakhanov_2009}. Optimal values of cost functional $J$, and  the number of iterations, line searches, and Newton steps in Table~\ref{Ex1_table1} are increasing for the higher values  of $V_{\max}$ because the linear system of equations to be solved become more stiff due to convection dominated character of the OCP problem.

\begin{figure}[t]
\centering
    \subfloat[Control $u$]{%
      \includegraphics[width=0.35\textwidth]{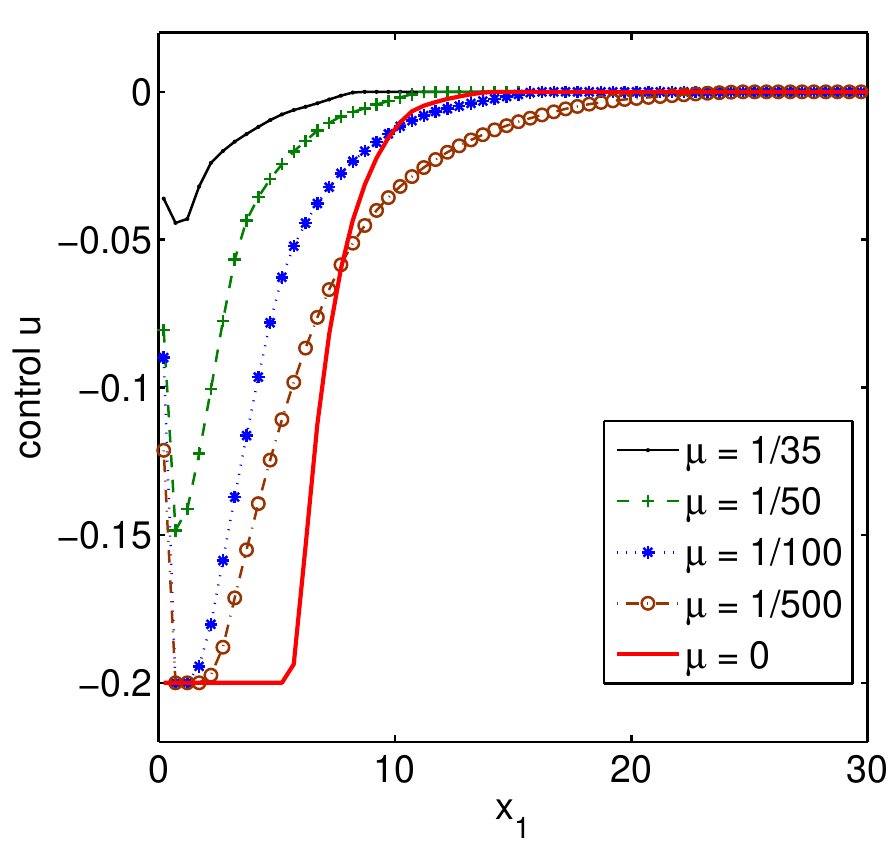}
    }
    \qquad \quad
    \subfloat[State $y$]{%
      \includegraphics[width=0.35\textwidth]{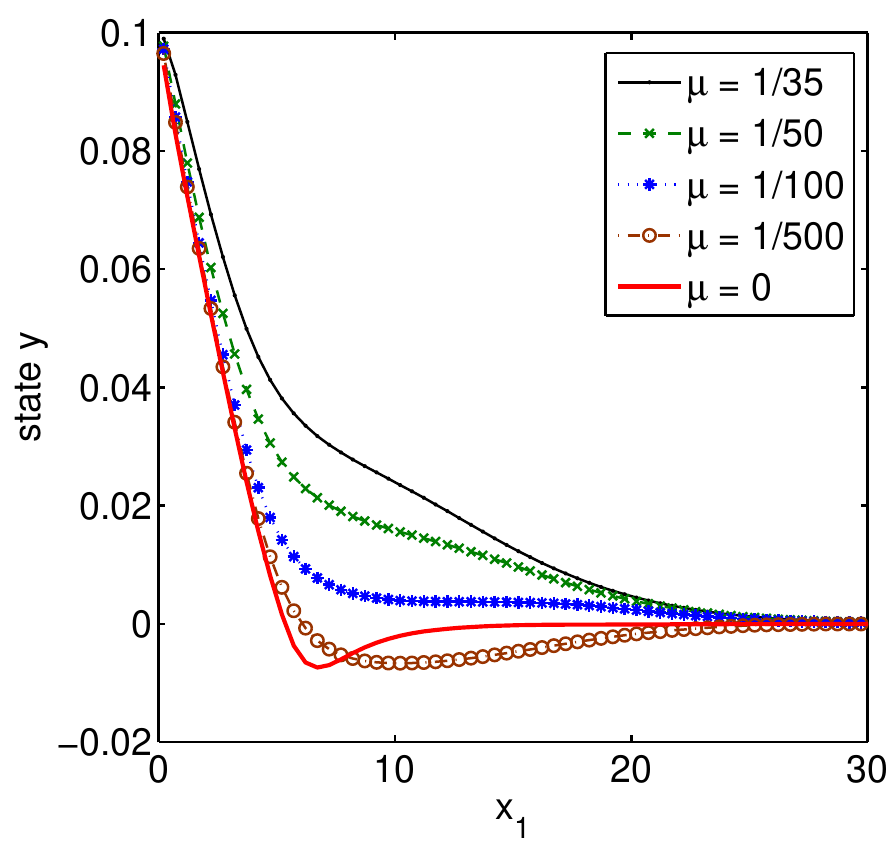}
    }
    \caption{Example~\ref{Ex1}: profiles of optimal controls $u$  and associated states $y$  at $t=0.75$ with $V_{\max}=32$ for various values of $\mu$.}
    \label{Fig:Ex1_sprP}
\end{figure}

\begin{table}[H]
\centering
  \begin{tabular}{r|c|c|r|r|r}
  $\mu$ & $J$ & $\mu j$ &$\#$ite.  & $\#$search  & $\#$Newton    \\
\hline
  $1/35$    & 1.10e-0 & 1.30e-1 &170    & 666  & 2193   \\
  $1/50$    & 7.31e-1 & 8.21e-2 &156    & 610  & 1971  \\
  $1/100$   & 3.60e-1 & 3.47e-2 &121     & 470  & 1516   \\
  $1/500$   & 1.21e-1 & 5.93e-3 &88     & 341  & 1099   \\
  $0$       & 5.80e-2 & 0     &77     & 297  & 956   \\
\hline
\end{tabular}
\caption{Example~\ref{Ex1}: optimal values of $J$, $\mu j$, and number of iterations, line searches, and Newton steps for $V_{max}=32$.}
\label{Ex1_table1s}
\end{table}

One of the important features of the optimal control is the robustness of the control with respect to the parameters. However, in the PDE-constrained optimization context, we do not always  have an explicit control function as in \cite{Peruzzi_2016,Tusset_2016}. In that case, we check the robustness of the control by solving the optimal control problem for different values of the parameters. For this example, we study the effect of the maximum velocity $V_{max}$ on the controlled wave solutions and of
the parameter $c_1$ on the stability of the waves. As $V_{max}$ increases, the  wave profiles are elongated $x_1$ directions and control is adjusted accordingly by shifting to the right in Figure~\ref{Fig:Ex1_const_Vmax}. Figure~\ref{Fig:Ex1_const_c1} shows that the solutions  display oscillations for higher values of $c_1$ as the uncontrolled  solutions in \cite{EAErmakova_EEShnol_MAPanteleev_AAButylin_VVolpert_FIAtaullakhanov_2009}. Since the control variable stays inside the bounds of the control, we  can conclude that the control variable is robust with respect to variations of the system parameters.

Now, we look into the effect of the sparse parameter $\mu$ on the optimization problem (\ref{M4}). Higher values of $\mu$ cause the sparsity of the optimal control, which can be clearly seen from  Figures~\ref{Fig:Ex1_spr} and \ref{Fig:Ex1_sprP}. Figure~\ref{Fig:Ex1_sprP} shows profiles of the controls $u$ and associated states $y$  along $x_1$ direction for various values of $\mu$. The values of the optimal cost functional and number of iterations with line searches and Newton steps are given in Table~\ref{Ex1_table1s}. When the sparse parameter $\mu$  increases, the smooth and non-smooth  cost functionals increase as well as the number nonlinear CG iterations, line searches and Newton iterations. However, we will show in the next example with terminal control, the number of  nonlinear CG iterations, line searches and Newton iterations can  decrease when the sparse parameter increases. Similar behavior was observed in \cite{ECasas_CRyll_FTroltzsch_2013} for different problems with sparse controls.

\begin{figure}[htp]
\centering
    \subfloat[$\mu=0$]{%
      \includegraphics[width=0.35\textwidth]{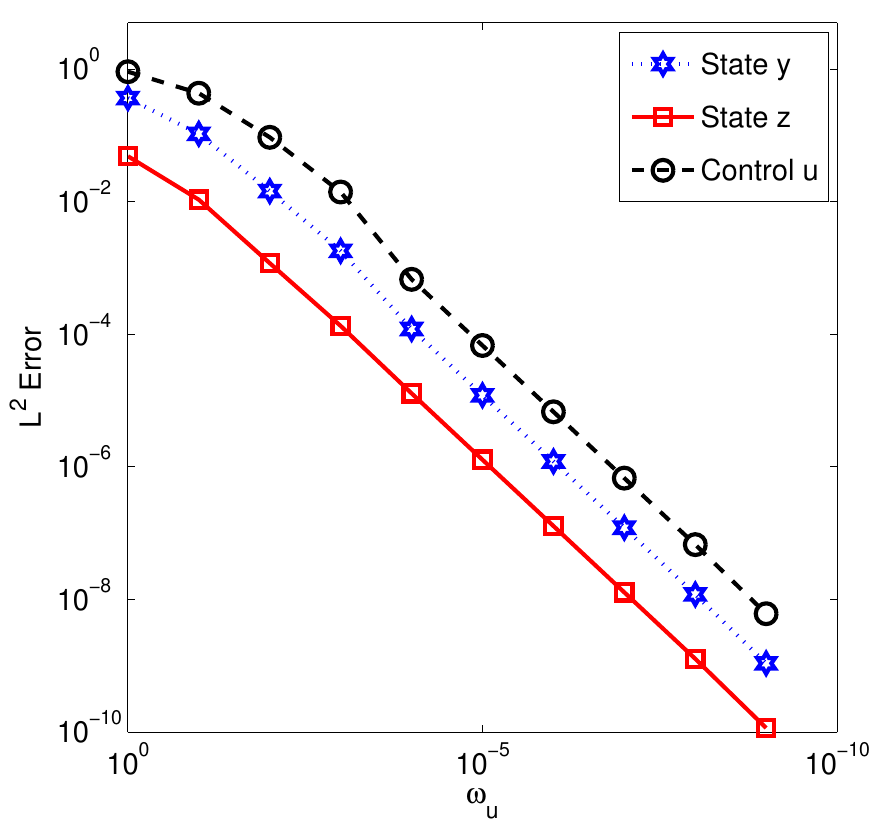}
    }
    \qquad
    \subfloat[$\mu=1/100$]{%
      \includegraphics[width=0.35\textwidth]{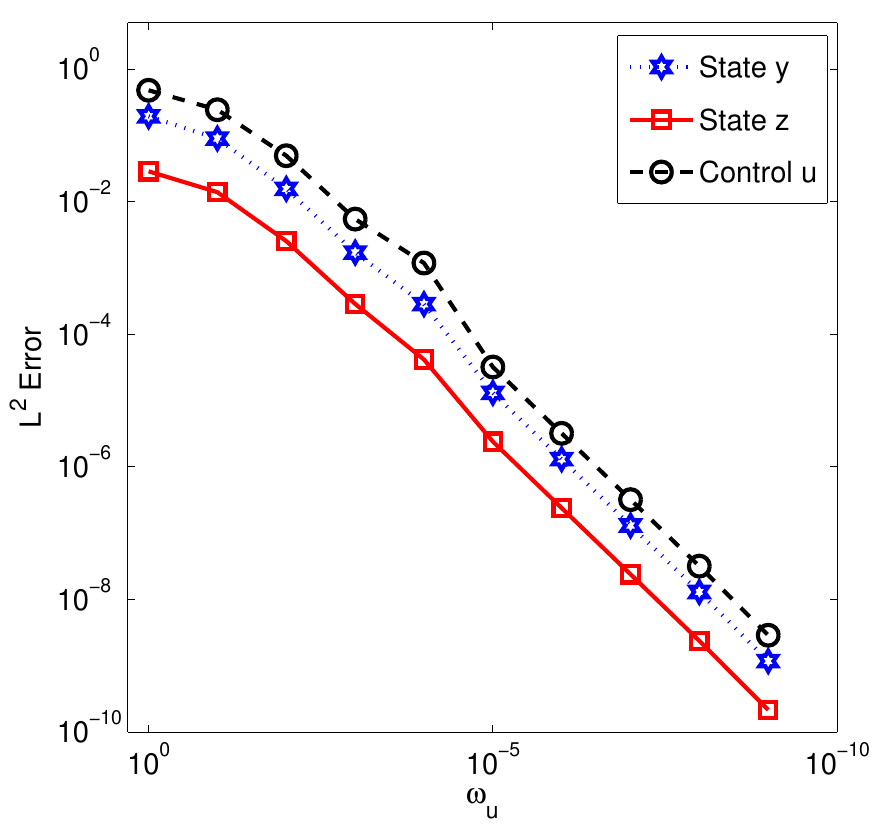}
    }
    \caption{Example~\ref{Ex1}: $L^2(Q)$  errors for $\|\overline{y}_{\omega_u}- \overline{y}_{\mathrm{ref}}\|$, $\|\overline{z}_{\omega_u}- \overline{z}_{\mathrm{ref}}\|$, and $\|\overline{u}_{\omega_u}- \overline{u}_{\mathrm{ref}}\|$.}
    \label{Fig:Ex1_regerror}
\end{figure}

Finally, in this example we study the convergence of the optimal control and its associated state for vanishing Tikhonov parameter  $\omega_u$ for  $\mu=0$ and $\mu \neq 0$. The convergence of optimal states and control for vanishing Tikhonov parameter  $\omega_u$  plays an important role for checking the SSCs according to the theory in
\cite{ECasas_CRyll_FTroltzsch_2015,ECasas_FTroltzsch_2016}. It was demonstrated for the one dimensional wave solutions of the classical FHN equation with sparse controls in   \cite{ECasas_CRyll_FTroltzsch_2015,CRyll_JLober_SMartens_HEngel_FTroltzsch_2016}. Here we obtain similar results  for two dimensional wave solutions of the convective FHN equations as shown in Figure~\ref{Fig:Ex1_regerror}.

We fix the maximum velocity $V_{\max}=64$. We take $\overline{y}_{\mathrm{ref}}:=y_{1e-10}$, $\overline{z}_{\mathrm{ref}}:=z_{1e-10}$, and $\overline{u}_{\mathrm{ref}}:=u_{1e-10}$ as reference solutions  to determine the order of convergence as
$\omega_u \downarrow 0$. We look for the errors, i.e., $\|\overline{y}_{\omega_u}- \overline{y}_{\mathrm{ref}}\|$, $\|\overline{z}_{\omega_u}- \overline{z}_{\mathrm{ref}}\|$, and $\|\overline{u}_{\omega_u}- \overline{u}_{\mathrm{ref}}\|$  as $\omega_u \downarrow 0$. We observe that the errors in $L^2$-norm decay linearly as $\omega_u \downarrow 0$ in  Figure~\ref{Fig:Ex1_regerror} for both cases $\mu=0$ and $\mu=1/100$. For the norm of $L^{\infty}(Q_T)$ this looks similar. The trajectories of the target functions, i.e., $y_Q$ and $z_Q$, are quite achieved as shown in Table~\ref{Ex1_table7}.

\begin{table}[t]
\centering
\begin{tabular}{c|c|c|c|c}
   &  \multicolumn{2}{c|}{$\mu=0$}   &  \multicolumn{2}{c}{$\mu=1/100$}  \\
   \cline{2-5}
  $\omega_u$ & $\|\overline{y}_{\omega_u}- y_Q \|_{L^2(Q)}$ & $\|\overline{z}_{\omega_u}- z_Q \|_{L^2(Q)}$ & $\|\overline{y}_{\omega_u}- y_Q \|_{L^{2}(Q)}$ & $\|\overline{z}_{\omega_c}- z_Q \|_{L^{2}(Q)}$ \\
\hline
  $1$      & 5.77167e-1 & 1.09495e-1 & 5.28730e-1 & 9.90557e-2  \\
  $1e-1$   & 3.84021e-1 & 7.32190e-2 & 4.48344e-1 & 8.48398e-2   \\
  $1e-2$   & 3.42134e-1 & 6.53493e-2 & 4.02207e-1 & 7.44464e-2  \\
  $1e-3$   & 3.39051e-1 & 6.46970e-2 & 3.95422e-1 & 7.25312e-2  \\
  $1e-4$   & 3.38605e-1 & 6.46413e-2 & 3.94826e-1 & 7.23321e-2  \\
  $1e-5$   & 3.38590e-1 & 6.46329e-2 & 3.94671e-1 & 7.22975e-2   \\
  $1e-6$   & 3.38588e-1 & 6.46321e-2 & 3.94668e-1 & 7.22959e-2  \\
  $1e-7$   & 3.38588e-1 & 6.46320e-2 & 3.94667e-1 & 7.22957e-2  \\
  $1e-8$   & 3.38588e-1 & 6.46320e-2 & 3.94667e-1 & 7.22957e-2  \\
  $1e-9$   & 3.38588e-1 & 6.46320e-2 & 3.94667e-1 & 7.22957e-2  \\
  $1e-10$  & 3.38588e-1 & 6.46320e-2 & 3.94667e-1 & 7.22957e-2  \\
\hline
\end{tabular}
\caption{Example~\ref{Ex1}: $L^2(Q)$  errors for $\|\overline{y}_{\omega_u}- y_Q \|$  and $\|\overline{z}_{\omega_u}- z_Q \|$.}
\label{Ex1_table7}
\end{table}

%%%%%%%%%%%%%%%%%%%%%%%%%%%%%%%%%%%%%%%%%%%%%%%%%%%%%%%%%%%%%%%%%%%%%%%%%%%
\subsection{Terminal control} \label{Ex2}

Our second test example is that the target functions, that is, $y_T(\cdot, T)$ and $z_T(\cdot,T)$, are only given at the final time. Regularization parameters are chosen as $\omega_{Q}^{y}=\omega_{Q}^{z}=0$, $\omega_{T}^{y}=\omega_{T}^{z}=1$,  and $\omega_u=10^{-3}$. We take the step sizes in space as $\Delta x_{1}=\Delta x_{2} =0.125$. We construct the desired functions as done in the previous example. They are given by
\[
y_T(x,T) = y_{\mathrm{nat}}(x,T/2) \quad \hbox{and} \quad z_T(x,T) = z_{\mathrm{nat}}(x,T/2),
\]
where $y_{\mathrm{nat}}$ and $z_{\mathrm{nat}}$ are the solutions of the uncontrolled convective FHN equation at the final time $T=1$ with the initial conditions
\[y_0(x,t) = \left\{
             \begin{array}{ll}
               1, & \hbox{if} \;\; x_a \leq x_1 \leq x_b,  \\
               0, & \hbox{otherwise},
             \end{array}
           \right.
\qquad
z_0(x,0) = 0,
\]
where $x_a =2$ and $x_b=2.2$. The admissible set of controls is defined as
\[
\mathcal{U}_{ad} := \{ u \in L^{\infty}(Q): \; 0 \leq u(x,t) \leq 0.2 \;\; \hbox{for   a.e   } (x,t) \in Q \}.
\]

\begin{figure}[H]
   \includegraphics[width=1\textwidth]{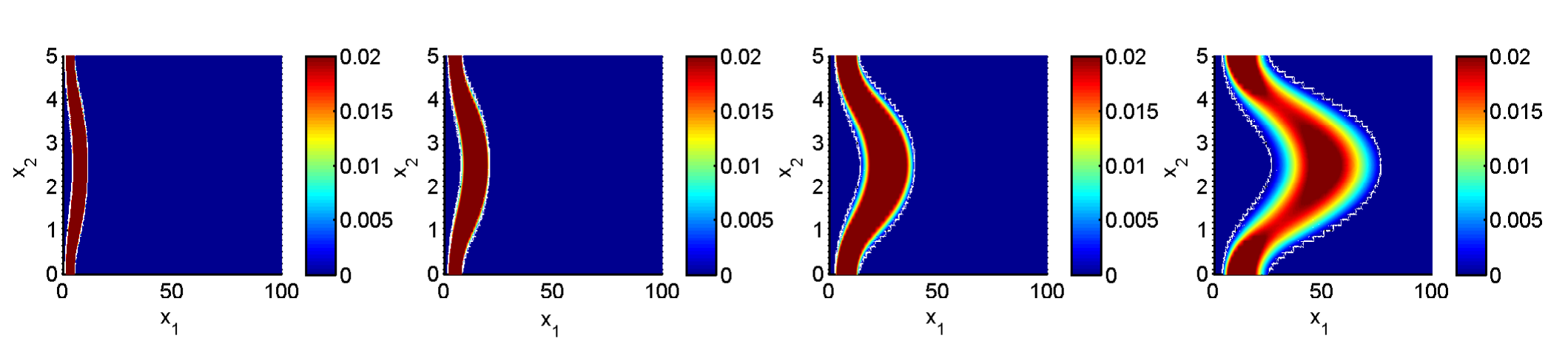}
   \includegraphics[width=1\textwidth]{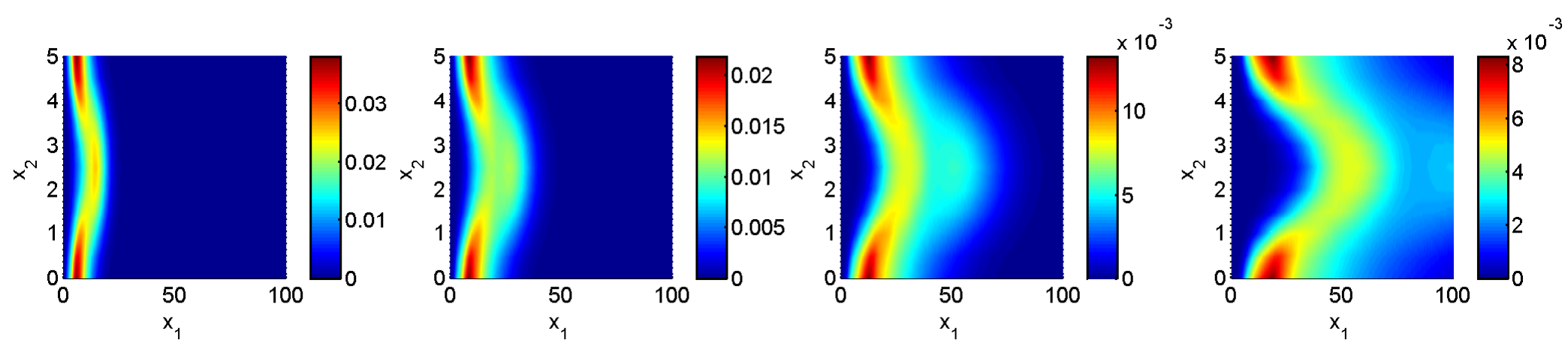}
    \caption{Example~\ref{Ex2}: optimal controls $u$ (top) and associated states $y$ (bottom) at $t=1$ without sparse control for  $V_{\max}$=16, 32, 64, 128 (from left to right).}
    \label{Fig:Ex2_const}
\end{figure}

\begin{table}[H]
\centering
\begin{tabular}{r|c|c|c|c}
  $V_{max}$ & $J$ &  $\#$ite.  & $\#$search  & $\#$Newton    \\
\hline
  $16$  & 1.12e-2 & 15    & 61  &  198   \\
  $32$  & 4.30e-2 & 16     & 65  & 211  \\
  $64$  & 1.42e-3 & 22     & 89  & 289   \\
  $128$  &2.89e-4 &25     & 95  & 304   \\
\hline
\end{tabular}
\caption{Example~\ref{Ex2}: optimal value of cost functional $J$, and  number of nonlinear CG iterations, line searches, and Newton steps without sparse control.}
\label{Ex2_table3}
\end{table}

We observe  similar behavior for the controls and associated states in Figure~\ref{Fig:Ex2_const}  and for the optimal value of cost functional $J$, the number of nonlinear CG iterations, line searches, and Newton steps  in Table~\ref{Ex2_table3}  as obtained for the whole space-time domain Example (\ref{Ex1}).

\begin{figure}[H]
   \includegraphics[width=1\textwidth]{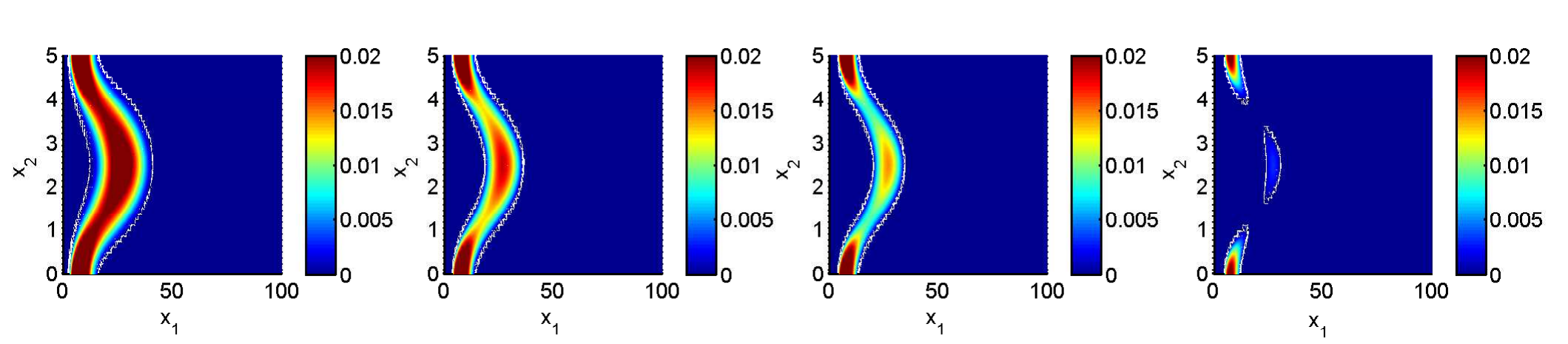}
   \includegraphics[width=1\textwidth]{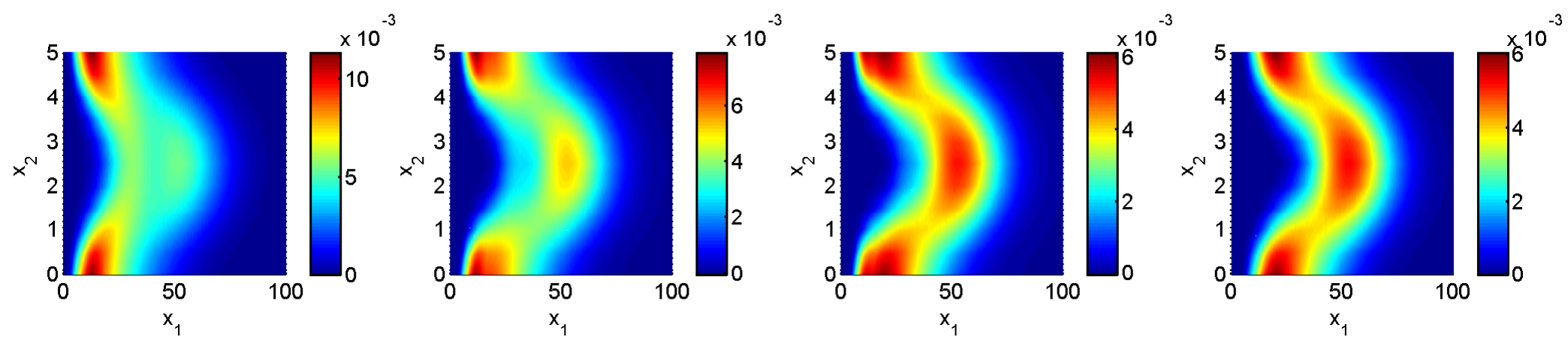}
    \caption{Example~\ref{Ex2}: optimal controls $u$ (top) and associated states $y$ (bottom) at $t=1$ for the sparse parameters $\mu=1/2000, 1/200, 1/150$, and $1/100$ (from left to right) and $V_{\max}=64$.}
    \label{Fig:Ex2_spr}
\end{figure}

\begin{figure}[H]
\centering
    \subfloat[Control $u$]{%
      \includegraphics[width=0.4\textwidth]{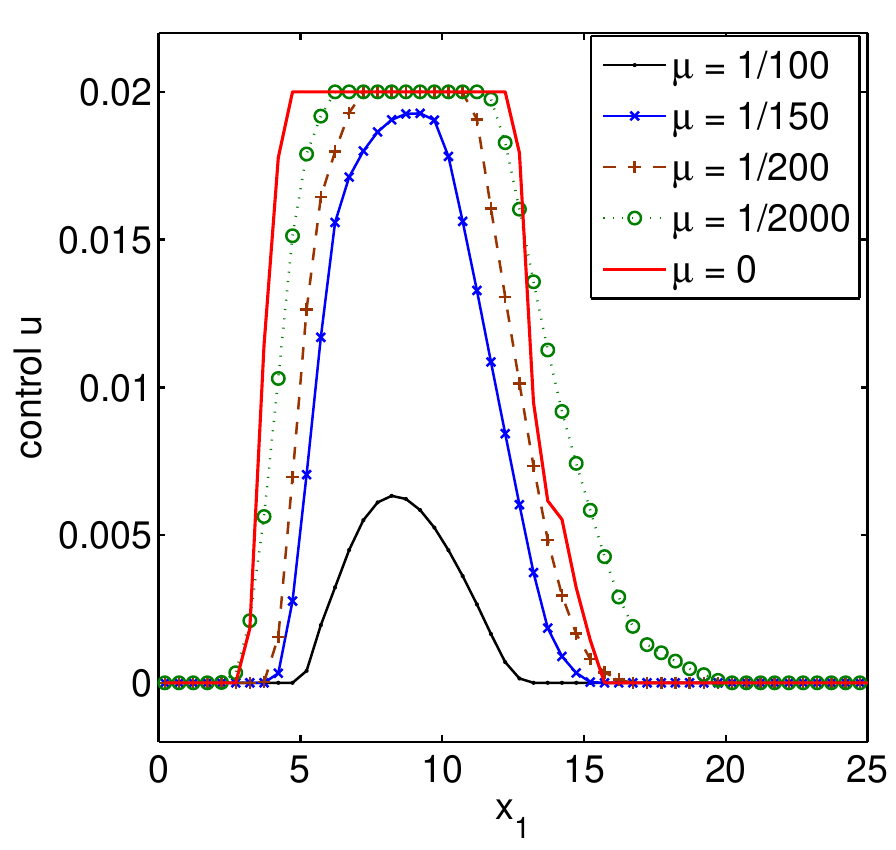}
    }
    \qquad \quad
    \subfloat[State $y$]{%
      \includegraphics[width=0.4\textwidth]{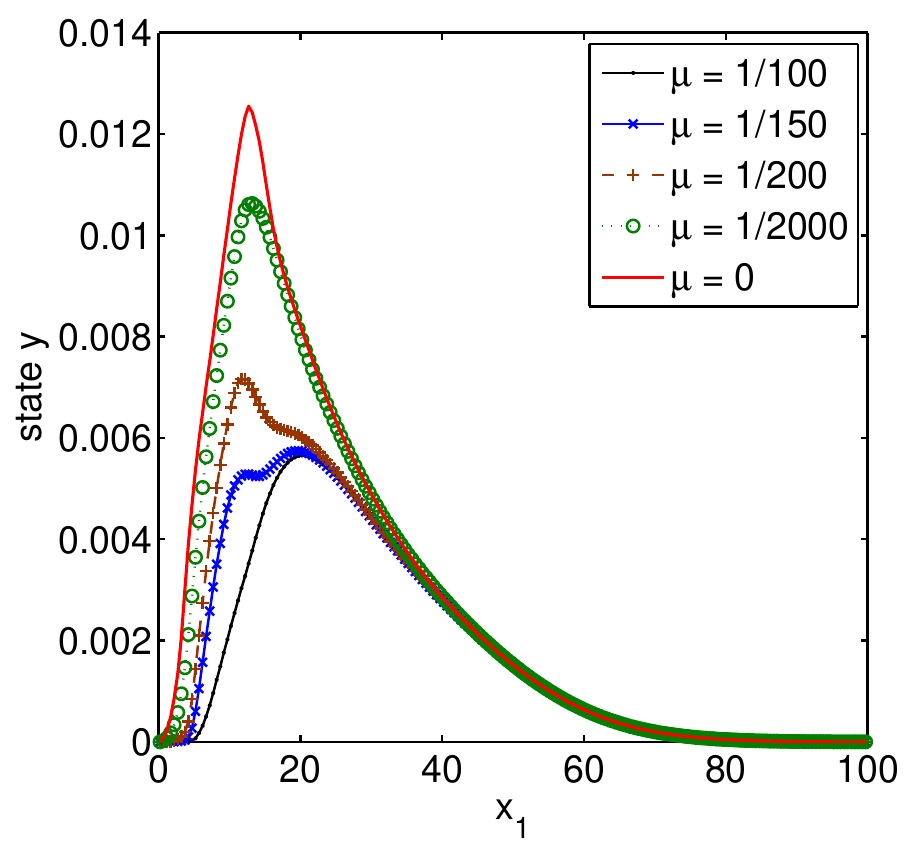}
    }
    \caption{Example~\ref{Ex2}: profiles of optimal controls $u$ and associated states $y$ (right) at $t=1$ for various values of the sparse parameter $\mu$.}
    \label{Fig:Ex2_sprT}
\end{figure}

\begin{figure}[H]
\centering
    \subfloat[$\mu=0$]{%
      \includegraphics[width=0.35\textwidth]{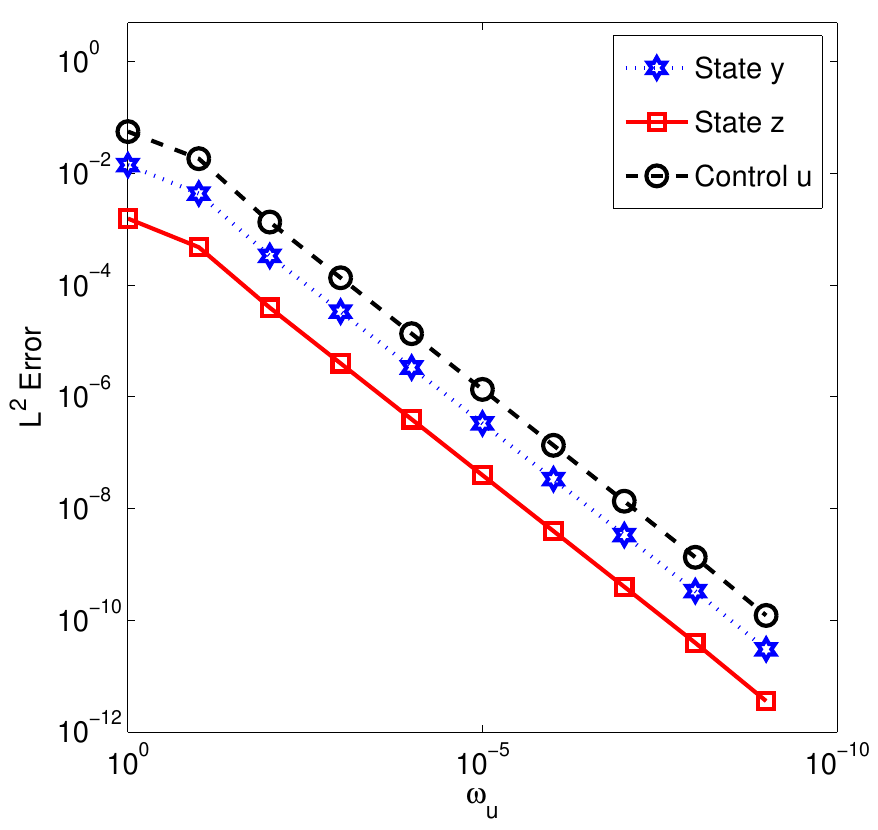}
    }
    \qquad
    \subfloat[$\mu=1/200$]{%
      \includegraphics[width=0.35\textwidth]{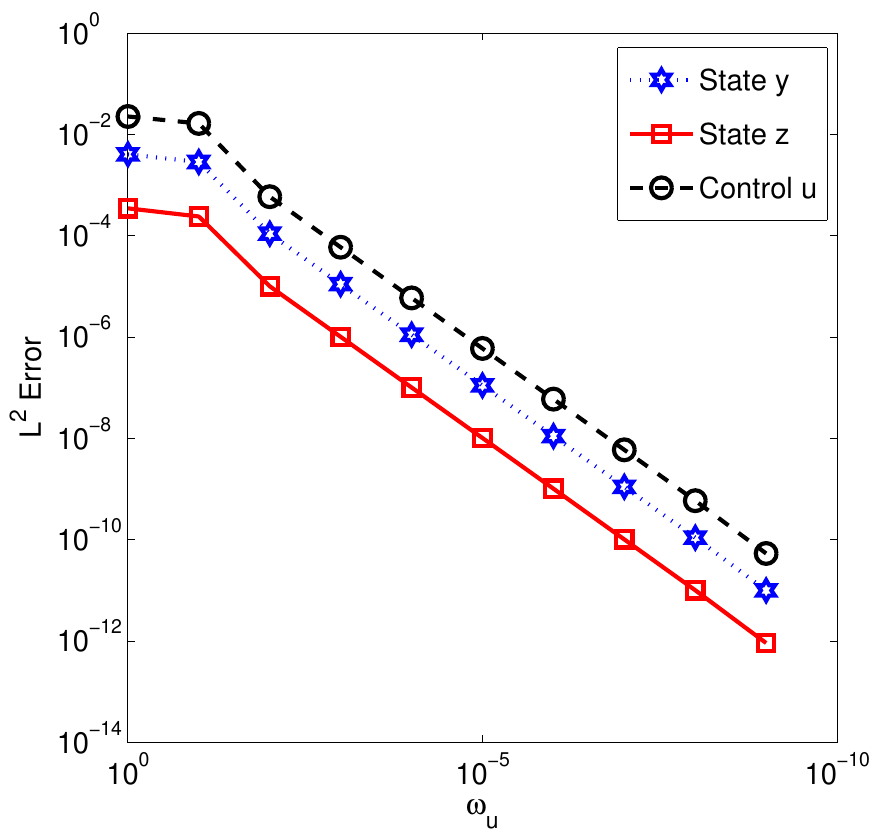}
    }
    \caption{Example~\ref{Ex2}: $L^2(Q)$  errors for $\|\overline{y}_{\omega_u}- \overline{y}_{\mathrm{ref}}\|$, $\|\overline{z}_{\omega_u}- \overline{z}_{\mathrm{ref}}\|$, and $\|\overline{u}_{\omega_u}- \overline{u}_{\mathrm{ref}}\|$.}
    \label{Fig:Ex2_regerror}
\end{figure}

The effect of the sparse parameter for the controlled solutions is displayed in Figure~\ref{Fig:Ex2_spr} and \ref{Fig:Ex2_sprT}. As stated before, in this case
for higher values of the sparse parameter $\mu$ the number of nonlinear CG iterations, line searches and Newton steps decrease as given in Table~\ref{Ex2_table5}.
We observe the same linear decay for the states and the control for vanishing Tikhonov parameter in Figure~\ref{Fig:Ex2_regerror} as in the whole space-time cylinder, cf. Figure~\ref{Fig:Ex1_regerror}. The values of  $\|\overline{y}_{\omega_u}- y_T \|$  and $\|\overline{z}_{\omega_u}- z_T \|$ in $L^2$-norm are given in Table~\ref{Ex2_table7}. We observe that the features of the desired trajectories are quite achieved.

Finally, in this example, we check the distance of discrete control $u_h$ from the local minima $u$ for the verification of  SSCs (\ref{ssc1}). Since we do not know the exact optimal solutions of this example, we take the discrete solutions  computed with $\Delta x_1 = \Delta x_2 = 0.3125$ as reference solutions. Table~\ref{Ex2_table8}  shows the numerical errors of $\|u-u_h\|$ and the error estimates given in (\ref{sscest}) for $\mu=1/200$. The numerical results verify the Theorem~\ref{Thm:est}. However, the error and estimator are not reduced sufficiently for finer discretizations.

\begin{table}[t]
\centering
\begin{tabular}{r|c|c|r|r|r}
  $\mu$ & $J_{opt}$ & $\mu j_{opt}$ &$\#$ite.  & $\#$search  & $\#$Newton    \\
\hline
  $1/100$   & 1.37e-1 & 1.15e-2 &2     & 5  & 16   \\
  $1/150$   & 9.24e-2 & 7.70e-3 &3     & 7  & 18   \\
  $1/200$   & 7.00e-2 & 5.88e-3 &4     & 8  & 19   \\
  $1/2000$  & 8.79e-3 & 6.43e-4 &19    & 74  & 238  \\
  $0$       & 1.42e-3 & 0       &22    & 89  & 289   \\
\hline
\end{tabular}
\caption{Example~\ref{Ex2}: optimal value of $J$, $\mu j$, and  number of nonlinear CG iterations, line searches, and Newton steps for $V_{\max}=64$.}
\label{Ex2_table5}
\end{table}

%%%%%%%%%%%%%%%%%%%%%%%%%%%%%%%%%%%%%%%%%%%%%%%%%%%%%%%%%%%%%%%%
\section{Conclusions}
Numerical results of the optimal control governed by the convective FHN model with traveling waves reveal different aspects of the parabolic semi-linear optimal control problems investigated. The second order optimality conditions for local solutions in form of 2D traveling waves are verified numerically for vanishing Tikhonov regularization parameter $\omega_u$  as done for one dimensional  waves of the classical FHN equation in \cite{ECasas_CRyll_FTroltzsch_2015,CRyll_JLober_SMartens_HEngel_FTroltzsch_2016}. By using the second order optimality conditions, we estimate the measure between the discrete solution and the local minima. The control of the convection dominated problems with wave solutions requires a large amount of computing time. In a future study we will investigate the reduced order modeling with proper orthogonal decomposition (POD) in space and model predictive control (MPC) in time \cite{Ryll14}.

\begin{table}[H]
\centering
\begin{tabular}{c|c|c|c|c}
   &  \multicolumn{2}{c|}{$\mu=0$}   &  \multicolumn{2}{c}{$\mu=1/200$}  \\
   \cline{2-5}
  $\omega_u$ & $\|\overline{y}_{\omega_u}- y_T \|_{L^2(\Omega)}$ & $\|\overline{z}_{\omega_u}- z_T \|_{L^2(\Omega)}$ & $\|\overline{y}_{\omega_u}- y_T \|_{L^{2}(\Omega)}$ & $\|\overline{z}_{\omega_c}- z_T \|_{L^{2}(\Omega)}$ \\
\hline
  $1$      & 6.53791e-2 & 2.77918e-2 & 7.92435e-2 & 2.81181e-2  \\
  $1e-1$   & 5.21169e-2 & 2.74520e-2 & 7.68586e-2 & 2.80154e-2   \\
  $1e-2$   & 4.87670e-2 & 2.73638e-2 & 7.13606e-2 & 2.77851e-2   \\
  $1e-3$   & 4.86308e-2 & 2.73684e-2 & 7.11792e-2 & 2.77772e-2   \\
  $1e-4$   & 4.86178e-2 & 2.73689e-2 & 7.11610e-2 & 2.77764e-2  \\
  $1e-5$   & 4.86165e-2 & 2.73690e-2 & 7.11591e-2 & 2.77763e-2  \\
  $1e-6$   & 4.86164e-2 & 2.73690e-2 & 7.11589e-2 & 2.77763e-2  \\
  $1e-7$   & 4.86164e-2 & 2.73690e-2 & 7.11589e-2 & 2.77763e-2  \\
  $1e-8$   & 4.86164e-2 & 2.73690e-2 & 7.11589e-2 & 2.77763e-2  \\
  $1e-9$   & 4.86164e-2 & 2.73690e-2 & 7.11589e-2 & 2.77763e-2  \\
  $1e-10$  & 4.86164e-2 & 2.73690e-2 & 7.11589e-2 & 2.77763e-2  \\
\hline
\end{tabular}
\caption{Example~\ref{Ex2}: $L^2(\Omega)$  errors for $\|\overline{y}_{\omega_u}- y_T \|$  and $\|\overline{z}_{\omega_u}- z_T \|$.}
\label{Ex2_table7}
\end{table}

\begin{table}[H]
\centering
\begin{tabular}{c|c|c}
  $\Delta x_1 = \Delta x_2$ & $\|u-u_h\|$ & $\frac{2}{\delta} \|\zeta\|$ \\
\hline
    $2.5$      & 2.43e-2 & 5.92e-2 \\
  $1.25$        & 2.43e-2 & 5.60e-2 \\
  $0.625$      & 2.46e-2 & 4.78e-2  \\
\hline
\end{tabular}
\caption{Example~\ref{Ex2}: Numerical errors of $\|u-u_h\|$ and error estimates $\frac{2}{\delta} \|\zeta\|$ for sparse controls with $\mu=1/200$.}
\label{Ex2_table8}
\end{table}

%We have used a discontinuous Galerkin method as a spatial discretization and a standard backward Euler method as a temporal discretization. Numerical test examples are presented fpr different set of parameters to control the travelling waves of the convective FHN sysem. We have also considered sparse optimal controls, which are applied only on a part of the spatial region and we tested the second order optimality condition for vanishing Tikhonov parameter.
%In spite of sparse controls, such optimal control problems can  require higher cpu time, especially for space-time target functions. Therefore, we address an implementation of model reduction techniques such as proper orthogonal diagonal decomposition with an adaptive strategy as a future work.

\section*{Acknowledgments}
We would like to thank Fredi Tr\"oltzsch for helpful discussions about this research. The authors also would like to express their sincere thanks to the referees for most valuable suggestions.

%\bibliographystyle{elsart-num-sort}
%\bibliographystyle{plain}
%\bibliography{references}

\end{document}